%% file: cas-sc-template.tex
\documentclass[a4paper,fleqn]{cas-sc}

\usepackage[authoryear,longnamesfirst]{natbib}
\usepackage{subcaption} 
\usepackage{siunitx}
\usepackage{ulem}

\def\tsc#1{\csdef{#1}{\textsc{\lowercase{#1}}\xspace}}
\tsc{WGM}
\tsc{QE}

\input{notation.tex}

\usepackage{xcolor}

\begin{document}
\let\WriteBookmarks\relax
\def\floatpagepagefraction{1}
\def\textpagefraction{.001}

\shorttitle{Wigner-Eckart Factorization of the Spectral Boltzmann Collision Operator}    

\shortauthors{R.R. Hiemstra, T. Ke\ss ler, and M.R.A. Abdelmalik}  

\title [mode = title]{Wigner-Eckart Factorization of the Spectral Boltzmann Collision Operator}  



\author[1,2]{René R. Hiemstra}[orcid=0000-0002-0179-9606]

\cormark[1]

\fnmark[1]

\ead{r.r.hiemstra@tue.nl}

\credit{Conceptualization, Methodology, Software, Validation, Writing}

\author[1,2]{Torsten Ke\ss ler}


\credit{Methodology, Review \& Editing}

\author[1,2]{Michael R.A. Abdelmalik}


\credit{Review \& Editing}

\affiliation[1]{organization={Department of Mechanical Engineering, Eindhoven University of Technology (TU/e)},
            country={The Netherlands}}
           
\affiliation[2]{organization={Simkinetic},
            url={www.simkinetic.tech},
            city={Veldhoven},
            country={The Netherlands}}

\cortext[1]{Corresponding author}

\fntext[1]{Marie Sk\l odowska-Curie Fellow.}


\begin{abstract}
We reduce the eight-dimensional weak form of the bilinear Boltzmann collision operator to a five-dimensional kinematic core by rigidly rotating the laboratory frame to align with the colliding pair and integrating over the $\mathrm{SO}(3)$ rotation group. This reduction yields an exact Wigner--Eckart factorization within a spectral Galerkin framework of associated Laguerre polynomials and spherical harmonics. The decomposition decouples the angular geometry from the scattering physics. The former, represented by Clebsch--Gordan coefficients, is evaluated exactly, while the latter is evaluated to machine precision by a spectrally convergent singular quadrature strategy. By explicitly zeroing specific entries, the macroscopic collision invariants are embedded without approximation. Cache-optimized contractions deliver up to a 37-fold single-core speedup and a 1000-fold memory reduction over standard dense Cartesian formulations. The approach is validated against analytical solutions for Maxwell molecules and infinite-order Chapman--Enskog viscosity coefficients for hard spheres.
\end{abstract}


\begin{highlights}
\item Wigner-Eckart factorization decouples scattering physics from angular geometry. 
\item 8D collision integral is exactly reduced to a 5D kinematic core via SO(3) symmetry.
\item Spectrally convergent numerical quadrature strategy for singular collision integral.
\item Null-space embedding natively conserves mass, momentum, and energy without approximation.
\item Compressed tensor storage reduces memory up to 99.9\% and accelerates runtime 37-fold.
\end{highlights}


\maketitle

\input{section1}
\input{section2}
\input{section3}
\input{section4}
\input{section5}
\input{section6}

\appendix
\input{appendix1}



\section*{Acknowledgements}
Funded by the European Union (Marie Skłodowska-Curie grant agreement No. 101105786). Views and opinions expressed are however those of the author(s) only and do not necessarily reflect those of the European Union or the European Research Executive Agency. Neither the European Union nor the granting authority can be held responsible for them.

\section*{Computer Code Availability}
The official software codebase is openly available on \href{https://github.com/simkinetic/collision-factorization}{GitHub}. The specific version used to generate the results in this paper is permanently archived on Zenodo under the DOI: \href{https://doi.org/10.5281/zenodo.20403094}{10.5281/zenodo.20403094}~\citep{hiemstra_zenodo_2026}.

\printcredits

\bibliographystyle{cas-model2-names}

\bibliography{cas-refs}

\end{document}

%% file: notation.tex
\newcommand{\kmax}{K_{\max}} 
\newcommand{\lmax}{L_{\max}} 
\newcommand{\state}{\alpha} 
\newcommand{\basis}{\phi} 
\newcommand{\radbasis}{R} 
\newcommand{\angbasis}{Y} 
\newcommand{\maxwellianref}{\tilde{M}} 

\newcommand{\eq}{\mathrm{eq}} 
\newcommand{\densityeq}{\rho_{\eq}} 
\newcommand{\densityref}{\tilde{\rho}} 
\newcommand{\tempeq}{T_{\eq}} 
\newcommand{\tempref}{\tilde{T}} 
\newcommand{\vmeaneq}{\mathbf{U}_{\eq}} 
\newcommand{\vmean}{\mathbf{U}} 

\newcommand{\vpre}{\mathbf{v}} 
\newcommand{\wpre}{\mathbf{w}} 
\newcommand{\vpost}{\mathbf{v}'} 
\newcommand{\wpost}{\mathbf{w}'} 

\newcommand{\vrelpre}{\mathbf{u}} 
\newcommand{\vrelpost}{\mathbf{u}'} 

\newcommand{\vpreunit}{\hat{\mathbf{v}}}    
\newcommand{\wpreunit}{\hat{\mathbf{w}}}    
\newcommand{\vpostunit}{\hat{\mathbf{v}}'}  
\newcommand{\scatterdir}{\hat{\mathbf{u}}'} 
\newcommand{\vreldir}{\hat{\mathbf{u}}} 

\newcommand{\incangle}{\beta}        
\newcommand{\defangle}{\chi}        
\newcommand{\aziangle}{\epsilon}    

\newcommand{\vmag}{v}               
\newcommand{\wmag}{w}               
\newcommand{\urel}{u}               

\newcommand{\Ctens}{C} 
\newcommand{\Rtens}{\mathcal{R}} 
\newcommand{\Gtens}{\mathcal{G}} 
\newcommand{\wignerD}{\mathcal{D}} 

\newcommand{\cmvar}{t} 
\newcommand{\srel}{s} 
\newcommand{\fracrel}{x} 
\newcommand{\cmenergy}{z} 
\newcommand{\halfangle}{y} 
\newcommand{\duffy}{h} 
\newcommand{\vhs}{\gamma} 
\newcommand{\chan}{c} 

\newcommand{\statepert}{\alpha_2}

\newcommand{\dof}{N} 
\newcommand{\numQ}{N_Q} 
\newcommand{\numK}{N_K} 
\newcommand{\numT}{N_T} 
\newcommand{\numG}{N_G} 
\newcommand{\numS}{N_S} 

%% file: section1.tex
\section{Introduction}
The Boltzmann equation provides the fundamental mesoscopic description of dilute gases far from equilibrium, where the collision operator $Q(f,f)$ governs the nonlinear evolution of the distribution function $f$. Accurately evaluating this five-dimensional integral remains the central computational bottleneck in deterministic kinetic simulations~\citep{pareschi:2000,gamba:2009,brechtken:2018}. State-of-the-art solvers generally employ either functional basis expansions or discrete velocity methods (DVMs). Fourier-based expansions~\citep{pareschi:2000,mouhot:2006,wu:2013,hu:2020} leverage fast convolutions but suffer from aliasing, domain truncation, and difficulties preserving macroscopic conservation laws (\cite{pareschi:2022}). Alternatively, spectral velocity-space expansions~\citep{gamba:2018,kessler:2019} may achieve spectral accuracy and exact conservation of invariant subspaces, but evaluating their dense Cartesian collision tensors incurs prohibitive $\mathcal{O}(L_{\max}^6)$ memory and arithmetic scaling with the maximum spherical harmonic degree, $L_{\max}$. Conversely, DVMs~\citep{cabannes:1980,platkowski:1988,bobylev:1995} avoid functional expansions entirely by directly discretizing the velocity space. While DVMs can be designed to preserve entropy and physical invariants~\citep{rogier:1994,brechtken:2018}, evaluating the multidimensional collision operator on a discrete grid introduces a steep tradeoff between accuracy and efficiency.

\subsection{Historical Context: Symmetry and Representation Theory}
The underlying difficulty of the standard Cartesian formulation is that it assembles the collision geometry in a fixed laboratory frame, thereby ignoring the intrinsic rotational invariance of binary elastic scattering. This tensorial structure has been recognized since the foundational works of \citet{chapman:1916} and \citet{enskog:1917}, and was later utilized by \citet{pekeris:1957} to extract the infinite-order limit of the Chapman--Enskog viscosity and heat-conduction coefficients for hard spheres.

Parallel developments in representation theory provided the framework to systematically decouple the collision integral. \citet{kumar:1967,kumar:1973} pioneered the use of irreducible spherical tensors and Racah--Wigner algebra to simplify the Chapman--Enskog method. Building on this tensorial foundation, \citet{massot:1976} applied the Wigner--Eckart theorem to prove the operator's reducibility across invariant subspaces. Early polar spectral schemes using associated Laguerre polynomials further highlighted the benefits of rotational symmetry~\citep{fonn:2014}, while invariance properties of Laguerre expansions were already noted by \citet{ender:1999}.

These symmetry-based insights remained largely analytical until the advent of modern spectral and Galerkin methods. The conservative Discontinuous Galerkin scheme of \citet{kitzler:2015, kitzler:2019}, for example, leveraged transformations to polar polynomial bases to reduce the evaluation complexity of the collision operator. Recent representation-theory approaches, such as those by \citet{hanke:2023} and \citet{cai:2020}, exploit $\mathrm{SO}(3)$ irreducible subspaces and solid spherical harmonics to achieve compression of the collision tensor, extending these solvers into transitional Knudsen regimes~\citep{wilkie:2023}. However, contemporary methods still face fundamental challenges in achieving accurate numerical quadrature, as well as overcoming the prohibitive memory and arithmetic costs associated with dense tensor contraction.

\subsection{Core Contributions}
In this work, we present a high-performance, conservative spectral method for the Boltzmann collision operator based on an exact dimensional reduction. Specifically, this paper introduces the following primary contributions:
\begin{itemize}

    \item \textbf{First-Principles Dimensional Reduction:} Unlike foundational algebraic approaches \citep{kumar:1967, massot:1976}, we derive the Wigner-Eckart factorization geometrically. By rigidly rotating the absolute laboratory frame into a privileged coordinate system aligned with the colliding pair and analytically integrating over the $\mathrm{SO}(3)$ rotation group, we decouple the 5-dimensional physical scattering tensor from the macroscopic angular geometry.
    
    \item \textbf{Spectrally Convergent Singular Quadrature:} To integrate the physical scattering tensor, we introduce a partitioning of the integration domain and coordinate transformation strategy that resolves the singular nature of the collision integral. This mapping yields a purely analytic integrand, allowing tailored Gauss-type quadrature rules to achieve spectral convergence and evaluate the entries of the tensor to machine precision.
    
    \item \textbf{Compressed Tensor Data Structure:} We develop a specialized, compressed data structure designed to exploit the Clebsch--Gordan selection rules governing the geometric Gaunt tensor. By storing the sparse macroscopic geometry in a flattened Coordinate (COO) format, we achieve memory savings ranging from 99\% to 99.9\% for increasing angular degrees, overcoming traditional storage bottlenecks.

    \item \textbf{Enforcement of Collision Invariants:} We identify a one-to-one mathematical relationship between specific discrete entries in the dense physical scattering tensor and the macroscopic collision invariants. By explicitly setting these radial transition entries to zero, we eliminate numerical quadrature errors, ensuring conservation of mass, momentum and energy.

    \item \textbf{Algorithmic Optimization and Validation:} We design and systematically analyze several tensor contraction algorithms tailored to the factorized data structure. Through benchmarking, we validate the implementation of the factorization and demonstrate that our cache-optimized contraction strategy yields single-core execution speedups of \(10\) to \(40\) relative to the standard dense formulation.
\end{itemize}

\subsection{Outline of the Paper}
Section~2 recalls the spatially homogeneous Boltzmann equation, the Galerkin method and the bilinear collision operator. Section~3 derives the five-dimensional reduction and the Wigner--Eckart factorization of the collision integral. Section~4 details the discrete implementation, including the singular quadrature strategy, the sparse storage of the Gaunt tensor, and cache-optimized contraction strategies. Section~5 presents numerical validation against the Bobylev--Krook--Wu solution, the Wang--Chang--Uhlenbeck eigenvalue spectrum, and the infinite-order Chapman--Enskog viscosity coefficients for hard-sphere molecules, alongside computational performance benchmarks. Finally, conclusions are drawn in Section~6.

%% file: section2.tex
\section{Background and Notation}
In this section, we establish the mathematical foundation and notation for the spatially homogeneous Boltzmann equation. We define the collision operator governing the particle interactions, introduce the spectral Galerkin method and formulate the discrete weak form of the collision tensor. Finally, we discuss the computational and storage bottleneck that arises from a standard Cartesian based description of the collision tensor, which motivates the dimensional reduction developed in Section 3.

Throughout this paper, we shall use boldface letters to denote vectors in $\mathbb{R}^3$ (e.g., $\vpre, \wpre$), their non-bold counterparts denote their scalar magnitudes or speeds (e.g., $\vmag = |\vpre|, \wmag = |\wpre|$), and hats denote unit directional vectors (e.g., $\vpreunit = \vpre/\vmag$).

\subsection{Model Problem}
We consider the initial value problem for the classical spatially homogeneous Boltzmann equation,
\begin{equation}\label{eq:BE}
	\frac{\partial f(t, \vpre)}{\partial t} = Q(f,f)(t, \vpre), \quad t \in \mathbb{R}^+, \quad \vpre \in \mathbb{R}^3,
\end{equation}
which describes the time evolution of the distribution function, $f : \mathbb{R}^+ \times \mathbb{R}^3 \to \mathbb{R}_{\geq 0}$, from its initial state $f(0, \vpre) = f_0(\vpre)$ to the final Maxwellian equilibrium distribution described by the conserved macroscopic density $\densityeq$, mean velocity $\vmeaneq$, and temperature $\tempeq$ of the gas.

The operator $Q(\cdot, \cdot)$ represents the rate of change of the distribution function due to binary elastic collisions between particles. For a colliding pair with pre-collision velocities $\vpre$ and $\wpre$, we define the relative velocity $\vrelpre = \vpre - \wpre$, possessing magnitude $\urel$ and direction $\vreldir$. Because elastic collisions conserve momentum and kinetic energy, the relative speed $\urel$ is preserved and the relative velocity rotates to a post-collision state $\vrelpost = \urel \scatterdir$, where the scattering direction $\scatterdir$ establishes the polar deflection angle $\cos \defangle = \vreldir \cdot \scatterdir$. The post-collision velocities are thereby uniquely determined by the center of mass velocity and the rotated relative velocity, yielding $\vpost = \frac{1}{2}(\vpre + \wpre + \vrelpost)$ and $\wpost = \frac{1}{2}(\vpre + \wpre - \vrelpost )$.

For a monatomic gas, the collision operator corresponds to the five-dimensional integral
\begin{equation}\label{eq:coll}
	Q(f,f)(\vpre) = \int_{\mathbb{R}^3} \int_{\mathbb{S}^2} B(\urel, \cos \defangle) \Bigl[ f(\vpost)f(\wpost) - f(\vpre)f(\wpre) \Bigr] \, d\scatterdir \, d\wpre.
\end{equation}
where the transition probability of the interaction is encapsulated within the physical collision kernel $B(\cdot, \cdot)$.

\subsection{Spectral Galerkin Bases}

To deterministically solve the spatially homogeneous Boltzmann equation, spectral methods expand the continuous distribution function onto a truncated basis of orthogonal polynomials. A natural choice for the three-dimensional velocity space is a basis composed of associated Laguerre polynomials \citep{chapman:1970} and spherical harmonics \citep{Varshalovich:1988}. This combination inherently captures the physical symmetries of the system: the spherical harmonics provide the rotational invariance of the Galerkin method, while the radial polynomials preserve the translational invariance of the collision operator. While the singular quadrature techniques developed later in this work rely specifically on the analytical properties of the Laguerre polynomials, the underlying factorization framework is geometric and applies to any velocity-space discretization that decouples into a generic radial basis function and spherical harmonics.

We define a composite spectral state index $\state = (k, l, m)$, where $k$ represents the radial polynomial degree, $l$ is the polar spherical harmonic degree, and $m$ is the azimuthal mode. The orthogonal basis functions are constructed by decoupling the radial and angular dependencies:
\begin{equation}
	\basis_\state(\vpre) = \radbasis_{k,l}(\vmag) \angbasis_{l,m}(\vpreunit).
\end{equation}
The distribution function is then expanded in terms of these basis functions and time-dependent spectral coefficients $c_\state(t)$:
\begin{equation}
	f(t, \vpre) \approx \maxwellianref(\vmag) \sum_{\state=1}^{N} c_\state(t) \basis_\state(\vpre),
\end{equation}	
where $N$ is the total number of basis functions in the functional expansion. The term $\maxwellianref(\vmag)$ is the reference Maxwellian weight function, defined by the reference density $\densityref$ and the reference temperature $\tempref$:
\begin{equation}
	\maxwellianref(\vmag) = \frac{\densityref}{(2\pi \tempref)^{3/2}} \exp\left(-\frac{\vmag^2}{2\tempref}\right).
\end{equation}

\subsection{Weak Form and Macroscopic Collision Tensor}
To obtain the weak form of the continuous Boltzmann equation, we multiply it by a test function $\basis_{\state_1}(\vpre)$ and integrate over the velocity space. Substituting the finite spectral expansion into this integral and applying the Galerkin approximation method transforms the integro-differential equation into a coupled system of ordinary differential equations. Because the orthogonal basis functions are normalized with respect to the Maxwellian weight $\maxwellianref(\vmag)$, the mass matrix on the left-hand side evaluates to the identity matrix, yielding the simplified system:
\begin{equation}
	\frac{\partial c_{\state_1}}{\partial t} = \sum_{\state_2=1}^{\dof} \sum_{\state_3=1}^{\dof} \Ctens_{\state_1, \state_2, \state_3} \, c_{\state_2} c_{\state_3}.
\end{equation}

The dynamics of this discrete system are governed by the macroscopic Cartesian collision tensor $\Ctens_{\state_1, \state_2, \state_3}$. To isolate the intrinsic scattering physics and angular geometry from the macroscopic flow parameters, we evaluate the collision tensor within a normalized reference frame where the characteristic density and temperature are scaled to unity ($\densityref = 1$, $\tempref = 1$). In this non-dimensionalized velocity space, the Maxwellian weight reduces to the normalized Gaussian distribution $\maxwellianref(\vmag) = (2\pi)^{-3/2} \exp(-\vmag^2/2)$. The physical dimensions are then recovered during the macroscopic time integration by appropriately scaling the collision operator.

Limiting our present scope to collision kernels that satisfy Grad's angular cut-off~\citep{grad:1958}, we can separate the continuous collision integral into its positive gain term and negative loss term. The total non-dimensional collision tensor is therefore expressed as $\Ctens_{\state_1, \state_2, \state_3} = \Ctens_{\state_1, \state_2, \state_3}^{+} - \Ctens_{\state_1, \state_2, \state_3}^{-}$. The elements of these respective tensors evaluate to eight-dimensional integrals over the absolute laboratory frame:
\begin{subequations}
\begin{align}
	\Ctens_{\state_1, \state_2, \state_3}^{+} &= \int_{\mathbb{R}^3} \int_{\mathbb{R}^3} \int_{\mathbb{S}^2} B(\urel, \cos \defangle) \basis_{\state_1}(\vpost) \basis_{\state_2}(\vpre) \basis_{\state_3}(\wpre) \maxwellianref(\vmag) \maxwellianref(\wmag) \, d\scatterdir \, d\wpre \, d\vpre, \\
	\Ctens_{\state_1, \state_2, \state_3}^{-} &= \int_{\mathbb{R}^3} \int_{\mathbb{R}^3} \int_{\mathbb{S}^2} B(\urel, \cos \defangle) \basis_{\state_1}(\vpre) \basis_{\state_2}(\vpre) \basis_{\state_3}(\wpre) \maxwellianref(\vmag) \maxwellianref(\wmag) \, d\scatterdir \, d\wpre \, d\vpre.
\end{align}
\end{subequations}

The evaluation of the macroscopic collision tensor $\Ctens$ represents a significant computational bottleneck in this spectral approach. If the basis is truncated at a maximum radial degree $\kmax$ and a maximum angular degree $\lmax$, the total number of degrees of freedom scales as $\dof \propto \kmax \lmax^2$. Consequently, the dense 3D Cartesian tensor contains $\dof^3 \propto \kmax^3 \lmax^6$ elements, which limits storage and scalability. However, this classical formulation evaluates the collision geometry within a fixed reference frame, thereby failing to exploit the rotational invariance of the binary scattering physics.

%% file: section3.tex
\section{Dimensional reduction of the collisional process}

This section details the reduction of the 8-dimensional Boltzmann collision integral to a 5-dimensional kinematic core. While previous representation-theory approaches achieved similar decouplings through abstract algebraic transformations~\citep{kumar:1967, massot:1976}, our formulation relies on a direct geometric approach. By integrating out the macroscopic spatial orientation over the $\mathrm{SO}(3)$ rotation group, the intrinsic collision process is decoupled from the macroscopic angular geometry, establishing the Wigner-Eckart factorization.

\subsection{Geometric Symmetries and Dimensionality Reduction}

\begin{figure}[htbp]
    \centering
    \includegraphics[width=0.45\textwidth]{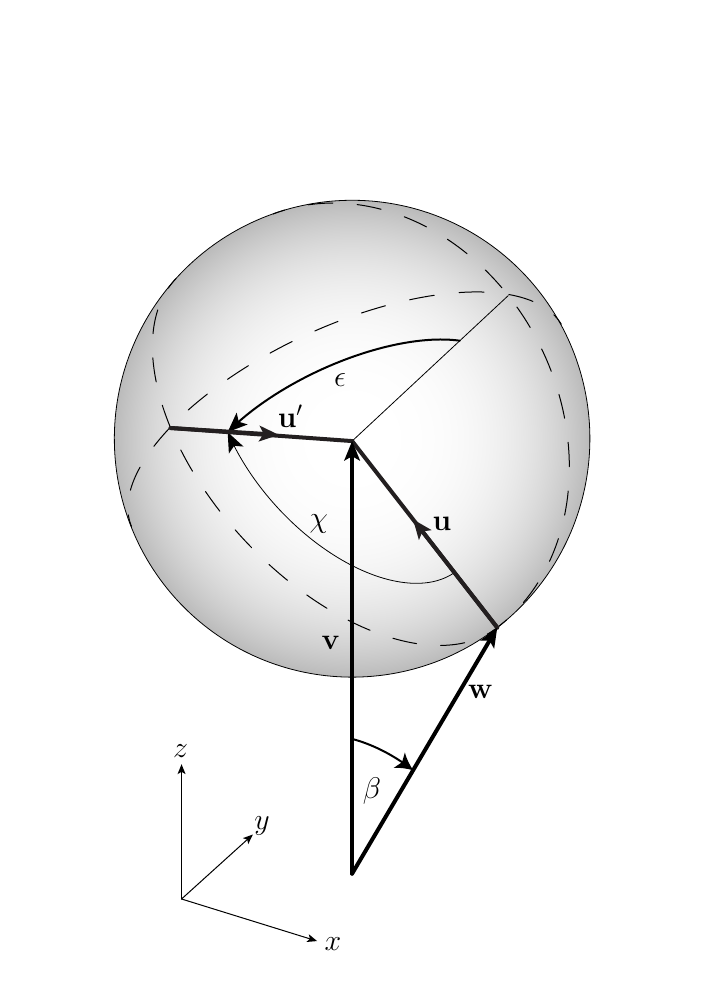}
    \caption{Geometric dimensional reduction of the 8-dimensional Cartesian collision integral to the privileged 5-dimensional kinematic core. The target velocity $\mathbf{v}$ is rigidly aligned with the local $z$-axis, and the incident velocity $\mathbf{w}$ is confined to the $x\text{-}z$ plane, separated by the polar incidence angle $\beta$. The pre-collision relative velocity $\mathbf{u}$ establishes the internal scattering axis. The post-collision relative velocity $\mathbf{u}'$ is parameterized on the scattering sphere by the polar deflection angle $\chi$ and the azimuthal angle $\epsilon$.}
    \label{fig:collision_frame}
\end{figure}

The evaluation of the Boltzmann collision operator requires integrating over an 8-dimensional space defined by $\vpre$, $\wpre$, and $\scatterdir$. However, the inherent isotropy of space dictates that the binary scattering kinematics possess no preferred spatial orientation. The collision process is governed by the scattering kernel $B(\urel, \cos \defangle)$, which depends on the relative speed $\urel$ and the polar deflection angle $\defangle$. Exploiting this rotational invariance eliminates these redundant degrees of freedom. We establish a privileged coordinate system by rigidly rotating the laboratory frame to align with the colliding pair (see Figure \ref{fig:collision_frame}):
\begin{enumerate}
    \item \textbf{Alignment of the target particle:} We define the local $z$-axis to be aligned with the target particle's velocity, such that $\vpre = \vmag\hat{\mathbf{z}}$. By locking $\vpre$ to the $z$-axis, we factor out the polar and azimuthal orientation of the collision event within the macroscopic reference frame. This eliminates two spatial degrees of freedom, reducing the integral from 8D to 6D.
    
    \item \textbf{Alignment of the incident particle:} The relative speed $\urel$ is invariant under rotations around the $z$-axis. Due to this cylindrical symmetry, we rotate our reference frame around $\vpre$ until the incident velocity $\wpre$ lies within the $x$-$z$ plane. This integrates out the azimuthal degree of freedom of $\wpre$, reducing the parameter space from 6D to 5D.
\end{enumerate}

With the coordinate system fixed, the geometric relationship between the pre- and post-collision states becomes constrained. Because both $\vpre$ and $\wpre$ lie within the $x$-$z$ plane, the pre-collision relative velocity $\vrelpre$ must also reside within this plane. 

During an elastic collision, conservation of momentum and kinetic energy dictates that the relative velocity rotates to a new post-collision state $\vrelpost$, while its magnitude is conserved ($|\vrelpost| = \urel$). Therefore, the vector $\vrelpost$ points to a location on a two-dimensional scattering sphere of radius $\urel$, centered at the origin of the relative frame.

Having factored out the three-dimensional rigid-body rotation of the collision pair, the intrinsic geometry of the elastic collision is described by five scalar variables:
\begin{itemize}
    \item $\vmag$: The speed of the target particle.
    \item $\wmag$: The speed of the incident particle.
    \item $\incangle$: The polar incidence angle between $\vpre$ and $\wpre$ (constrained to the $x$-$z$ plane).
    \item $\defangle$: The polar deflection angle between the pre- and post-collision relative velocities ($\vrelpre$ and $\vrelpost$).
    \item $\aziangle$: The azimuthal scattering angle of $\vrelpost$ rotating around the axis defined by $\vrelpre$.
\end{itemize}

Any approximation of the Boltzmann operator using a spectral basis must ultimately evaluate the collision dynamics over these five intrinsic dimensions.

\subsection{Exact Factorization via Integration over the Rotation Group}
Assuming the angular basis functions $\angbasis_{l m}$ are complex spherical harmonics, any collision configuration in the laboratory frame is generated by applying a 3D rigid rotation $R \in \mathrm{SO}(3)$ to the privileged 5D configuration established above. 

Expressing the 8-dimensional Cartesian measure in spherical coordinates yields,
\begin{equation}
	d\vpre \, d\wpre \, d\scatterdir = (\vmag^2 \, d\vmag \, d\vpreunit) (\wmag^2 \, d\wmag \, d\wpreunit) d\scatterdir.
\end{equation}
The spatial orientation of the unit vector pair $(\vpreunit, \wpreunit)$ requires four angular coordinates. Let $\vpreunit \in S^2$ be defined in the global frame by its polar and azimuthal angles $(\theta, \phi)$, yielding the standard solid angle $d\vpreunit = \sin\theta \, d\theta \, d\phi$. By treating $\vpreunit$ as the $z$-axis of a local reference frame, $\wpreunit$ is parameterized by the internal polar incidence angle $\incangle$ and a local azimuthal angle $\psi$, giving $d\wpreunit = \sin\incangle \, d\incangle \, d\psi$.

The triplet $(\phi, \theta, \psi)$ defines a standard Z-Y-Z Euler rotation $R \in \mathrm{SO}(3)$ that governs the absolute rigid-body orientation of the collision pair. By rearranging the direct product of the solid angles, we isolate the invariant Haar measure $dR = \sin\theta \, d\theta \, d\phi \, d\psi$, which yields the geometric relation,
\begin{equation}
	d\vpreunit \, d\wpreunit = dR \sin\incangle \, d\incangle.
\end{equation}
Because the scalar speeds $(\vmag, \wmag)$ and the local scattering solid angle $(d\scatterdir = \sin\defangle \, d\defangle \, d\aziangle)$ are invariant under the global rotation $R$, the integration measure mathematically decouples the invariant internal scattering geometry from the arbitrary reference frame,
\begin{equation}
	d\vpre \, d\wpre \, d\scatterdir = dR \, d\Phi_{5D}, \quad \text{with} \quad d\Phi_{5D} = \vmag^2 \wmag^2 \sin\incangle \sin\defangle \, d\vmag \, d\wmag \, d\incangle \, d\defangle \, d\aziangle,
\end{equation}
where $dR$ is the standard Haar measure over the rotation group $\mathrm{SO}(3)$, normalized such that $\int dR = 8\pi^2$~\citep{Varshalovich:1988}, and $d\Phi_{5D}$ encompasses the 5-dimensional kinematic core.

Under a rigid rotation $R$, the radial polynomials $\radbasis_{k_i l_i}$ and the physical scattering kernel $B(\urel, \cos\defangle)$ remain invariant. However, the complex spherical harmonics transform linearly according to the Wigner D-matrices $\wignerD^{l}_{m m'}(R)$. We can therefore express the laboratory spherical harmonics in terms of the rotation $R$ acting on the local privileged frame directions (denoted with a subscript 0),
\begin{equation}
    \angbasis_{lm}(\hat{\mathbf{n}}) = \sum_{m'} \wignerD^{l}_{m m'}(R) \angbasis_{lm'}(\hat{\mathbf{n}}_0).
\end{equation}

Substituting these transformed harmonics into the collision integral decouples the integration over the global rigid-body rotation ($dR$) from the invariant 5D kinematic core ($d\Phi_{5D}$). Isolating the gain term yields,
\begin{equation}
    \Ctens_{\state_1, \state_2, \state_3}^{+} = \sum_{m_1', m_2', m_3'} \underbrace{\left( \int_{\mathrm{SO}(3)} \wignerD^{l_1}_{m_1 m_1'}(R) \wignerD^{l_2}_{m_2 m_2'}(R) \wignerD^{l_3}_{m_3 m_3'}(R) dR \right)}_{\text{Macroscopic Geometry}} \times \underbrace{\int_{5D} (\quad \dots \quad) \, d\Phi_{5D}}_{\text{Intrinsic Collision Process}}.
\end{equation}

To evaluate the macroscopic geometry integral, we apply the fundamental representation identity for the integration of three Wigner D-matrices over the rotation group~\citep{Varshalovich:1988},
\begin{equation}
    \int_{\mathrm{SO}(3)} \wignerD^{l_1}_{m_1 m_1'}(R) \wignerD^{l_2}_{m_2 m_2'}(R) \wignerD^{l_3}_{m_3 m_3'}(R) dR = 8\pi^2 \begin{pmatrix} l_1 & l_2 & l_3 \\ m_1 & m_2 & m_3 \end{pmatrix} \begin{pmatrix} l_1 & l_2 & l_3 \\ m_1' & m_2' & m_3' \end{pmatrix}.
\end{equation}

This mathematical identity is the origin of the Wigner-Eckart factorization. It reveals that the macroscopic rotational geometry splits into two distinct angular couplings: a global Wigner 3-$j$ symbol depending on the laboratory modes ($m_i$), and a local Wigner 3-$j$ symbol coupling to the intrinsic kinematic modes ($m_i'$). Because the Wigner 3-$j$ symbols possess closed-form analytical representations in terms of finite sums of factorials~\citep{Varshalovich:1988}, this macroscopic angular geometry can be evaluated without approximation.

The global 3-$j$ symbol isolates the macroscopic angular geometry, while the local 3-$j$ symbol is absorbed into the 5-dimensional integral to weight the physical scattering dynamics. Consequently, the gain term $\Ctens^+$ separates into a sparse macroscopic angular component $\Gtens$ and a dense kinematic tensor $\Rtens^+$,
\begin{equation}
    \Ctens_{\state_1, \state_2, \state_3}^{+} = \Gtens_{m_1 m_2 m_3}^{l_1 l_2 l_3} \, \Rtens_{k_1 k_2 k_3}^{+l_1 l_2 l_3}.
\end{equation}
Here, $\Gtens$ is the macroscopic Gaunt tensor, defined directly by the global Wigner 3-$j$ symbol,
\begin{equation}
    \Gtens_{m_1 m_2 m_3}^{l_1 l_2 l_3} = 8\pi^2 \begin{pmatrix} l_1 & l_2 & l_3 \\ m_1 & m_2 & m_3 \end{pmatrix},
\end{equation}
and $\Rtens$ is the reduced physical tensor, which pulls the remaining summation over the local modes into the 5-dimensional kinematic integral,
\begin{equation} \label{eq:R_tensor_full}
    \Rtens_{k_1 k_2 k_3}^{+l_1 l_2 l_3} = \int_{5D} B(\urel, \cos\defangle) \radbasis_{k_1 l_1}(\vmag') \radbasis_{k_2 l_2}(\vmag) \radbasis_{k_3 l_3}(\wmag) P_{\text{gain}}^{l_1 l_2 l_3}(\vpostunit_0, \vpreunit_0, \wpreunit_0) \, \maxwellianref(\vmag) \maxwellianref(\wmag) \, d\Phi_{5D}.
\end{equation}

The term $P_{\text{gain}}^{l_1 l_2 l_3}$ represents the invariant geometric filter, aggregating the local spherical harmonics weighted by the surviving local 3-$j$ symbol,
\begin{equation}
    P_{\text{gain}}^{l_1 l_2 l_3}(\vpostunit_0, \vpreunit_0, \wpreunit_0) = \sum_{m_1', m_2', m_3'} \begin{pmatrix} l_1 & l_2 & l_3 \\ m_1' & m_2' & m_3' \end{pmatrix} \angbasis_{l_1 m_1'}(\vpostunit_0) \angbasis_{l_2 m_2'}(\vpreunit_0) \angbasis_{l_3 m_3'}(\wpreunit_0).
\end{equation}

Furthermore, the geometric constraints of our privileged frame simplify this internal summation. Because the target particle velocity is locked to the local $z$-axis ($\vpreunit_0 = \hat{\mathbf{z}}$), its spherical harmonic evaluates to zero for all modes except $m_2' = 0$. Azimuthal conservation within the local 3-$j$ symbol enforces $m_3' = -m_1'$. Letting $m' = m_1'$, the gain filter evaluates as,
\begin{equation}
	P_{\text{gain}}^{l_1 l_2 l_3}(\vpostunit_0, \incangle) = \angbasis_{l_2 0}(\hat{\mathbf{z}}) \sum_{m'} \begin{pmatrix} l_1 & l_2 & l_3 \\ m' & 0 & -m' \end{pmatrix} \angbasis_{l_1 m'}(\vpostunit_0) \angbasis_{l_3, -m'}(\incangle, 0).
\end{equation}

\subsection{Incorporation of the Loss Term}
While the preceding derivation details the positive gain term $\Ctens^{+}$, the negative loss term $\Ctens^{-}$ follows the same rotational decoupling. For the loss term, however, the test function evaluates the target particle at its pre-collision state. In our privileged frame, this state is already locked to the local $z$-axis, which restricts its spherical harmonic to $m_1' = 0$. 

Returning to the generic invariant geometric filter, the spherical harmonics for both the test function and the target particle now evaluate to zero for all modes except $m_1' = 0$ and $m_2' = 0$. Azimuthal conservation within the local 3-$j$ symbol ($m_1' + m_2' + m_3' = 0$) forces the incident particle's mode to be zero ($m_3' = 0$). Consequently, the internal summation collapses. The geometric filter for the loss term reduces to a zonal mode:
\begin{equation}
	P_{\text{loss}}^{l_1 l_2 l_3}(\incangle) = \angbasis_{l_1 0}(\hat{\mathbf{z}}) \angbasis_{l_2 0}(\hat{\mathbf{z}}) \begin{pmatrix} l_1 & l_2 & l_3 \\ 0 & 0 & 0 \end{pmatrix} \angbasis_{l_3 0}(\incangle, 0).
\end{equation}
This formulation demonstrates why the loss kinematics are independent of the scattering angles $\defangle$ and $\aziangle$; the only remaining angular dependence is the polar incidence angle $\incangle$.

In practice, we do not evaluate or store the loss term as a separate operator. Instead, the loss kinematics are combined with the gain kinematics under the same 5-dimensional integral. The total reduced physical tensor $\Rtens = \Rtens^{+} - \Rtens^{-}$ is computed directly:
\begin{equation}
	\Rtens_{k_1 k_2 k_3}^{l_1 l_2 l_3} = \int_{5D} B(\urel, \cos\defangle) \radbasis_{k_2 l_2}(\vmag) \radbasis_{k_3 l_3}(\wmag) \left[ \radbasis_{k_1 l_1}(\vmag') P_{\text{gain}}^{l_1 l_2 l_3}(\vpostunit_0, \incangle) - \radbasis_{k_1 l_1}(\vmag) P_{\text{loss}}^{l_1 l_2 l_3}(\incangle) \right] \maxwellianref(\vmag) \maxwellianref(\wmag) \, d\Phi_{5D}.
\end{equation}
By evaluating the difference between the gain and loss terms directly inside the 5D kinematic core, the full nonlinear collision operator is compressed into a single dense physical tensor $\Rtens$.

%% file: section4.tex
\section{Numerical Implementation and Algorithm Design}
This section translates the continuous Wigner-Eckart factorization into a discrete computational framework and is available as an open-source Matlab library with C++ extension \citep{hiemstra_zenodo_2026}. The numerical implementation is structured into three programmatic phases. First, we discretize the macroscopic angular geometry, exploiting geometric selection rules to build a compressed, sparse routing table. Second, we compute the dense 5D collision tensor, introducing a singular quadrature strategy that resolves the non-analytic variable hard-sphere scattering process while maintaining spectral convergence, followed by a null-space embedding to guarantee macroscopic conservation. Finally, we detail cache-optimized tensor contraction strategies that fuse these bifurcated data structures, overcoming traditional memory and arithmetic bottlenecks.

\subsection{Discrete Index Mapping and Real-Valued Transformations}
To translate the continuous Wigner-Eckart factorization into a numerical framework, we must map the continuous spectral modes to one-dimensional contiguous memory indices. We decompose the composite state $\state_i = (k_i, l_i, m_i)$ into a radial dimension and a flattened angular dimension.

First, we define the composite angular index $q_i$, which maps a specific spherical harmonic degree and mode $(l_i, m_i)$ to a 1D integer array. For a given maximum angular resolution $\lmax$, the index $q_i$ enumerates all valid spherical harmonics, running from $1$ to $\numQ = (\lmax+1)^2$. 

The macroscopic geometry captured by the Gaunt tensor $\Gtens$ represents the integration of three spherical harmonics over the unit sphere. Allocating a dense 3D grid for all possible polar degrees $(l_1, l_2, l_3)$ is inefficient. Instead, evaluating this geometry in the complex Spherical Harmonics basis yields Wigner 3-$j$ symbols. The complex (denoted by subscript $\mathbb{C}$) geometric coupling $\Gtens_{\mathbb{C}}$ evaluates to zero unless all of the following geometric selection rules are simultaneously satisfied~\citep{Varshalovich:1988}:
\begin{enumerate}
    \item \textbf{Triangle Inequality:} The polar degrees must form a closed vector triangle: $|l_2 - l_3| \le l_1 \le l_2 + l_3$.
    \item \textbf{Parity Constraint:} The sum of the polar degrees $l_1 + l_2 + l_3$ must be an even integer.
    \item \textbf{Azimuthal Balance:} The azimuthal modes satisfy: $m_1 = m_2 + m_3$.
\end{enumerate}
Because these rules forbid the majority of interactions, we extract only the valid triplets $(l_1, l_2, l_3)$ and assign them a 1D interaction channel index $\chan \in \{1, \dots, \numT\}$. 

While the complex basis provides clean selection rules, simulating physical gas flows requires real-valued distribution functions. We therefore map the complex tensor to the Real Spherical Harmonic ($\mathbb{R}$) basis via a unitary transformation matrix $U$~\citep{ivanic:1996}. The macroscopic real tensor is obtained via the tensor contraction,
\begin{equation}
    \Gtens_{\mathbb{R}} = (U \otimes U \otimes U^*) \Gtens_{\mathbb{C}}.
\end{equation}
Because the unitary matrix $U$ mixes conjugate azimuthal modes ($+m$ and $-m$) within the same polar degree $l$, the transformation alters the specific weights of the Gaunt coefficients but leaves the interaction channel $\chan = (l_1, l_2, l_3)$ invariant. This yields a real-valued discrete factorization,
\begin{equation}
    \Ctens_{q_1, q_2, q_3}^{k_1, k_2, k_3} = \Rtens_{k_1, k_2, k_3}^{(\chan)} \times [\Gtens_{\mathbb{R}}]_{q_1, q_2, q_3}^{(\chan)}.
    \label{eq:wigner-eckart-factorization}
\end{equation}

With the mathematical decoupling established and the interaction channels $\chan$ identified, the factorization naturally bifurcates the operator into two complementary memory structures. The physical tensor \(\Rtens\) depends only on the bounded radial limits (\(k_i \le \kmax\)) and the dense valid channels \(\chan\). It is therefore constructed as a dense 4D contiguous array, \(\Rtens[k_1, k_2, k_3, \chan]\). Conversely, the macroscopic geometry $\Gtens_{\mathbb{R}}$ is sparse and requires a specialized routing architecture, which we detail next.

\subsection{Coordinate (COO) Memory Layout}
To exploit the sparsity of the macroscopic geometry $\Gtens_{\mathbb{R}}$, we avoid allocating a dense 3D grid for the angular states and instead store the tensor in a flattened Coordinate (COO) format. We define a linear memory index $z \in \{1, \dots, \numG\}$ that iterates over the non-zero geometric collisions. 

As illustrated in Table \ref{tab:coo_layout}, each row $z$ in the COO structure acts as a discrete routing instruction. It identifies the three interacting macroscopic angular states $(q_1, q_2, q_3)$, provides the scalar geometric weight $g_z$, and supplies the channel index $\chan$. 

\begin{table}[h]
\centering
\caption{Sparse coordinate (COO) memory layout of the real Gaunt tensor $\Gtens_{\mathbb{R}}$ for a macroscopic angular resolution of $\lmax = 4$. Out of $\numQ^3 = 25^3 = 15,625$ possible macroscopic angular triplets, only $\numG = 1158$ satisfy the geometric selection rules. The index $z$ iterates strictly over these non-zero elements, acting as a sparse routing table that links the macroscopic angular states $q_i$ to the appropriate physical interaction channel $\chan \in \{1, \dots, 42\}$.}
\label{tab:coo_layout}
\begin{tabular}{@{}cccccc@{}}
\toprule
\textbf{Memory Offset ($z$)} & \textbf{Test State ($q_1$)} & \textbf{Target ($q_2$)} & \textbf{Incident ($q_3$)} & \textbf{Channel ($\chan$)} & \textbf{Real Gaunt Weight ($g_z$)} \\ \midrule
1                            & 1                           & 1                       & 1                         & 1                      & 0.2821                             \\
2                            & 2                           & 2                       & 1                         & 7                      & 0.2821                             \\
3                            & 3                           & 3                       & 1                         & 7                      & 0.2821                             \\
4                            & 4                           & 4                       & 1                         & 7                      & 0.2821                             \\
$\vdots$                     & $\vdots$                    & $\vdots$                & $\vdots$                  & $\vdots$               & $\vdots$                           \\
28                           & 7                           & 2                       & 2                         & 15                     & -0.1262                            \\
$\vdots$                     & $\vdots$                    & $\vdots$                & $\vdots$                  & $\vdots$               & $\vdots$                           \\
1158                         & 21                          & 25                      & 25                        & 42                     & 0.1065                             \\ \bottomrule
\end{tabular}
\end{table}

\subsection{Singular Quadrature Architecture} \label{sec:ingularquad}
While the macroscopic geometry is sparsely routed via the COO list, the intrinsic collision physics must be evaluated to populate the dense 4D physical array $\Rtens^{(\chan)}$. The Wigner-Eckart factorization reduces this physical collision tensor to a five-dimensional kinematic core governed by the intrinsic variables $\vmag$, $\wmag$, $\incangle$, $\defangle$, and $\aziangle$. Evaluating this reduced tensor requires integrating a decoupled geometric filter against the physical scattering dynamics. However, evaluating this integral directly in standard speeds introduces a function-space mismatch due to the non-analytic behavior of a typical cross section ($\urel^{\vhs}$) across the $\vmag=\wmag$ diagonal, preventing the realization of spectral convergence. To resolve this singularity, we decompose the domain of integration and introduce specialized coordinate transformations that align with the geometry of the singular sub-domain. 

\subsubsection{Rotated Kinematic Decoupling and the Cone Singularity}
To decouple the kinematics and absorb the singularities, we rigidly rotate the radial domain by 45 degrees, substituting the original speeds with center-of-mass and relative coordinates: $\cmvar = (\vmag+\wmag)/\sqrt{2}$ and $\srel = (\vmag-\wmag)/\sqrt{2}$. We restrict the domain to the upper half-plane ($\vmag>\wmag$) and map these coordinates to spectral integration domains via an energy variable $\cmenergy = \cmvar^2 \in [0, \infty)$ and a fractional relative coordinate $\fracrel = \srel/\cmvar \in [0, 1]$. Concurrently, the polar incidence angle is parameterized via a non-linear half-angle transformation $\halfangle = \sin(\incangle/2) \in [0, 1]$.

While this transformation factors the base Maxwellian envelope $\maxwellianref$ without cross-terms and eliminates the infinite first derivative at $\incangle \to 0$, the relative speed magnitude $\urel = \sqrt{2\cmenergy}\sqrt{\fracrel^2 + (1-\fracrel^2)\halfangle^2}$ introduces a new 2D cone singularity, $\sqrt{\fracrel^2+\halfangle^2}$, at the origin $(\fracrel,\halfangle) \to (0,0)$. 

\subsubsection{The 2D Duffy Transformation}
To resolve the cone singularity, the coupled variables must be algebraically factored out of the square root. We apply a 2D Duffy transformation, splitting the $(\fracrel, \halfangle)$ unit square along the diagonal into two triangular sub-domains mapped by a local coordinate $\duffy \in [0, 1]$:

\begin{itemize}
    \item \textbf{Patch 1: Lower Triangle ($\fracrel > \halfangle$):} We map $\halfangle = \fracrel \cdot \duffy$. The relative speed analytically factors the $\fracrel$ variable out of the singularity, yielding $\urel = \sqrt{2\cmenergy} \cdot \fracrel \sqrt{1 + (1-\fracrel^2)\duffy^2}$.
    \item \textbf{Patch 2: Upper Triangle ($\halfangle > \fracrel$):} We map $\fracrel = \halfangle \cdot \duffy$. The relative speed factors the $\halfangle$ variable out, yielding $\urel = \sqrt{2\cmenergy} \cdot \halfangle \sqrt{1 + (1-\halfangle^2) \duffy^2}$.
\end{itemize}

By precomputing these factored limits, the numerical square root never evaluates a sum near zero. This mapping regularizes the infinite derivatives while preserving orthogonal $[0,1]^2$ integration boundaries and the tensor-product topology of the quadrature grid.

\subsubsection{5D Tensor-Product Quadrature Grid and Spectral Convergence}
By combining the rotated kinematic decoupling with the Duffy transformation, the full 5D kinematic core is integrated using a multi-dimensional tensor product rule across the two triangular patches. The grid for each patch is constructed from five specific Gauss-type rules, each chosen for its native weight function to resolve the underlying physical and geometric behavior of its respective coordinate,

\begin{itemize}
    \item \textbf{Center-of-Mass Energy ($\cmenergy \in [0,\infty)$):} Integrated using a \textbf{Generalized Gauss-Laguerre} rule with the weight function $W(\cmenergy) = \cmenergy^{\vhs/2}e^{-\cmenergy}$. By setting the parameter $\vhs$ to match the singularity order of the scattering kernel, the weight explicitly accounts for the continuous algebraic singularity. The physical singularity ($\urel^\vhs \propto \cmenergy^{\vhs/2}$) is canceled by the integration weight, while the background Maxwellian's exponential decay ensures the evaluated remainder is bounded globally.
    \item \textbf{Primary Fractional Coordinate ($\fracrel$ or $\halfangle \in [0, 1]$):} Depending on the active Duffy patch, the outer integration spans either the fractional relative speed $\fracrel$ or the incidence half-angle $\halfangle$. Integrated using a \textbf{Gauss-Legendre} rule. This integrates the smooth polynomial interactions of the Cartesian measure and the solid spherical harmonic basis functions.
    \item \textbf{Local Duffy Coordinate ($\duffy \in [0, 1]$):} Integrated using a \textbf{Gauss-Legendre} rule. This evaluates the internal mapping coordinate that factors out the 2D cone singularity. The residual square roots in the factored relative speed $\urel$ are evaluated away from zero, ensuring the integrand remains analytic.
    \item \textbf{Polar Deflection Angle ($\cos\defangle \in [-1, 1]$):} Integrated using a \textbf{Gauss-Legendre} rule. This accurately resolves the smooth polynomial scattering dependencies of the physical collision kernel over the finite interval.
    \item \textbf{Azimuthal Scattering Angle ($\aziangle \in [0, 2\pi]$):} Integrated using a \textbf{Trapezoidal Rule}. The trapezoidal rule functions as the exact Gauss-type rule for periodic domains, providing spectrally exact integration for the periodic trigonometric functions spanning the azimuthal rotation.
\end{itemize}

To maintain the spectral convergence of the Galerkin framework, the quadrature grid must resolve the algebraic degree dictated by the state truncation limits $\kmax$ and $\lmax$. The minimum number of points required to integrate the polynomial components of the transformed Cartesian measure is detailed in Appendix \ref{app:quad_bounds} and is summarized in Table \ref{tab:quadrature_bounds}. 

\begin{table}[ht]
\centering
\caption{Minimum theoretical quadrature bounds required to achieve exact spectral integration of the polynomial components within the 5D kinematic core. The resolutions are dictated by the radial ($\kmax$) and angular ($\lmax$) spectral truncation limits.}
\label{tab:quadrature_bounds}
\begin{tabular}{@{}llllc@{}}
\toprule
\textbf{Coordinate} & \textbf{Active Domain} & \textbf{Quadrature Rule} & \textbf{Weight Function} & \textbf{Minimum Points ($N$)} \\ \midrule
Energy ($\cmenergy$) & Both Patches & Gen. Gauss-Laguerre & $W(\cmenergy) = \cmenergy^{\vhs/2}e^{-\cmenergy}$ & $\lceil \frac{3\kmax + 1.5\lmax + 3}{2} \rceil$ \\
Outer Frac. ($\fracrel$) & Patch 1 ($\fracrel > \halfangle$) & Gauss-Legendre & $W(\fracrel) = 1$ & $4\kmax + 3\lmax + 4$ \\
Outer Angle ($\halfangle$) & Patch 2 ($\halfangle > \fracrel$) & Gauss-Legendre & $W(\halfangle) = 1$ & $4\kmax + 3\lmax + 4$ \\
Inner Duffy ($\duffy$) & Patch 1 ($\halfangle = \fracrel \duffy$) & Gauss-Legendre & $W(\duffy) = 1$ & $\kmax + \lceil 1.5\lmax \rceil + 1$ \\
Inner Duffy ($\duffy$) & Patch 2 ($\fracrel = \halfangle \duffy$) & Gauss-Legendre & $W(\duffy) = 1$ & $3\kmax + \lfloor 1.5\lmax \rfloor + 3$ \\
Deflection ($\cos\defangle$) & Both Patches & Gauss-Legendre & $W(\cos\defangle) = 1$ & $\kmax + \lceil 0.5\lmax \rceil + 1$ \\
Azimuthal ($\aziangle$) & Both Patches & Trapezoidal & $W(\aziangle) = 1$ & $2\kmax + \lmax + 1$ \\ \bottomrule
\end{tabular}
\end{table}

Furthermore, we achieve machine precision for the non-polynomial components (the residual Gaussian envelope and the fractional scattering kernel) by augmenting these baseline polynomial bounds with targeted padding (10 - 20 additional points). Because the rotational scattering angles ($\defangle, \aziangle$) remain isolated from the non-polynomial singularities, this padding is applied to the coupled kinematic coordinates ($\cmenergy, \fracrel, \halfangle, \duffy$).

Finally, the discrete evaluation over these domains is optimized via a block-decoupled map-reduce architecture and sum factorization. An exposition of this is presented in greater detail in Appendix \ref{app:1}.

\subsection{Enforcement of Macroscopic Conservation}
To guarantee conservation, we exploit a mathematical mapping between specific discrete entries in the dense physical tensor $\Rtens$ and the macroscopic collision invariants.

\subsubsection{Mapping Invariants to Spectral Indices}
In the Boltzmann equation, the collision integral evaluates to zero when tested against the collision invariants: $\psi(\vpre) \in \{1, \vpre, |\vpre|^2\}$. In the chosen Laguerre-spherical-harmonic basis, these physical invariants map directly to the lowest-order functions,
\begin{itemize}
    \item \textbf{Mass ($1$):} The density invariant is an isotropic constant. In the spectral basis, this corresponds to the radial degree $k_1 = 0$ and the polar degree $l_1 = 0$.
    \item \textbf{Momentum ($\vpre$):} The linear velocity term represents a dipole. This corresponds to the radial degree $k_1 = 0$ and the polar degree $l_1 = 1$, spanning the three azimuthal modes $m_1 \in \{-1, 0, 1\}$.
    \item \textbf{Energy ($|\vpre|^2$):} The kinetic energy invariant is an isotropic quadratic term. In the Laguerre basis, this corresponds to the radial degree $k_1 = 1$ and the polar degree $l_1 = 0$.
\end{itemize}

\subsubsection{The Null-Space Correction}
To correct the discrete operator against these quadrature errors, we iterate through the precomputed physical tensor $\Rtens$ and overwrite the integrated values with exact zeros for the following specific slices:

\begin{enumerate}
    \item \textbf{Mass Conservation:} For all interaction channels $\chan$ where $l_1 = 0$, we set:
    \begin{equation}
        \Rtens_{0, k_2, k_3}^{(\chan)} = 0 \quad \forall \ k_2, k_3
    \end{equation}
    \item \textbf{Momentum Conservation:} For all interaction channels $\chan$ where $l_1 = 1$, we set:
    \begin{equation}
        \Rtens_{0, k_2, k_3}^{(\chan)} = 0 \quad \forall \ k_2, k_3
    \end{equation}
    \item \textbf{Energy Conservation:} For all interaction channels $\chan$ where $l_1 = 0$, we set:
    \begin{equation}
        \Rtens_{1, k_2, k_3}^{(\chan)} = 0 \quad \forall \ k_2, k_3
    \end{equation}
\end{enumerate}

Because the Gaunt tensor $\Gtens_{\mathbb{R}}$ acts as a routing table, setting these specific dense blocks in $\Rtens$ to zero guarantees that any geometric transition pointing to the macroscopic states $(k_1=0, l_1=0)$, $(k_1=0, l_1=1)$, or $(k_1=1, l_1=0)$ will be multiplied by zero during the runtime contraction. The physical null space is thereby "baked in" to the compressed data structure itself, ensuring machine-precision conservation during the temporal evolution of the gas without incurring any runtime computational overhead.

\subsubsection{Enforcement of Detailed Balance}
We enforce detailed balance at equilibrium by zeroing out the corresponding input states in the precomputed physical tensor. For the specific interaction channel $\chan$ where $l_2 = 0$ and $l_3 = 0$ (which forces $l_1 = 0$), we set the resulting vector to zero:
\begin{equation}
    \Rtens_{k_1, 0, 0}^{(\chan)} = 0 \quad \forall k_1.
\end{equation}
By zeroing these entries we guarantee that our background $\maxwellianref$ satisfies $Q (\maxwellianref, \maxwellianref ) = 0$.

\subsection{Tensor Contraction Strategies}
With both the sparse geometric routing table $\Gtens_{\mathbb{R}}$ constructed and the dense physical tensor $\Rtens$ integrated and corrected, the discrete operator is assembled. The optimal evaluation strategy for executing this factorized collision operator during a simulation is contingent on the target hardware architecture (e.g., shared-memory CPUs vs.\ massively parallel GPUs) and the specific use case. For instance, evaluating the linearized collision operator allows for pre-computation of the absolute equilibrium state, while spatial batch-processing (computing the operator simultaneously across thousands of physical grid points) is advantageous on GPU architectures thanks to coalesced memory access and massive parallelism. The present algorithmic analysis focuses on the foundational baseline: evaluating the full, nonlinear bilinear collision operator at a single spatial point on a standard shared-memory CPU. 

To evaluate the nonlinear collision integral $Q_{k_1, q_1}$ for a specific test state, the continuous representation is replaced by the discrete factored summation:
\begin{equation}
    Q_{k_1, q_1} = \sum_{z=1}^{\numG} g_z \sum_{k_2=0}^{\kmax} \sum_{k_3=0}^{\kmax} \Rtens_{k_1, k_2, k_3}^{(\chan_z)} f_{k_2, q_{2z}} f_{k_3, q_{3z}} \label{eq:full_discrete_sum}
\end{equation}
where $z$ iterates over the non-zero rows of the sparse geometric COO list, fetching the geometric weight $g_z$, the interacting angular states $(q_{2z}, q_{3z})$, and the corresponding physical channel $\chan_z$. Because addition is commutative, the summation loops in Equation \ref{eq:full_discrete_sum} can be algebraically rearranged. In addition to the dense contraction (baseline), we detail three distinct tensor contraction strategies, each exhibiting a profoundly different hardware memory profile.

\subsubsection{Dense Cartesian Contraction (Baseline)}
The standard spectral approach abandons the internal physical symmetries of the collision process and evaluates the operator using the fully assembled 3D Cartesian tensor $\Ctens_{\state_1, \state_2, \state_3}$. The nonlinear collision integral is evaluated via a dense summation over the spectral basis,
\begin{equation}
    Q_{\state_1} = \sum_{\state_2=1}^{\dof} \sum_{\state_3=1}^{\dof} \Ctens_{\state_1, \state_2, \state_3} f_{\state_2} f_{\state_3}
\end{equation}
where $\dof = (\kmax+1)(\lmax+1)^2$ is the total number of spectral degrees of freedom. While this formulation provides contiguous memory reads and a trivial algebraic implementation, the operation scales as $\mathcal{O}(\dof^3) \sim \mathcal{O}(\kmax^3 \lmax^6)$. Consequently, this approach is presented as the theoretical baseline against which the efficiency of the Wigner-Eckart factorized strategies will be measured.

\subsubsection{The Naive Contraction}
The naive approach evaluates Equation \ref{eq:full_discrete_sum} as written. The algorithm linearly streams through the 1D geometric array $z \in \{1, \dots, \numG\}$. For every valid geometric collision, it fetches the appropriate channel block $\Rtens^{(\chan_z)}$ and executes the dense $K^3$ radial energy summation. 

This approach couples the sparse memory routing directly to the dense floating-point arithmetic. Because the indices $q_{2z}$ and $q_{3z}$ pseudo-randomly jump across the spectral velocity space as $z$ increments, the algorithm repeatedly incurs cache misses when fetching the distribution coefficients $f_{k, q}$, limiting CPU throughput.

\subsubsection{Radial-First Contraction (Pairing $\Rtens$ first)}
Conversely, one can elect to compute the physical energy exchange for all possible trial pairs prior to routing them geometrically. In this approach, we contract the dense physical tensor $\Rtens$ with the distribution functions to build a dense intermediate state tensor $\Psi$,
\begin{equation}
    \Psi_{q_2, q_3, \chan}^{k_1} = \sum_{k_2, k_3} \Rtens_{k_1, k_2, k_3}^{(\chan)} f_{k_2, q_2} f_{k_3, q_3}
\end{equation}
Following this dense pre-computation, the final collision operator is obtained via a single, geometric streaming pass over the COO list,
\begin{equation}
    Q_{k_1, q_1} = \sum_{\{z \mid q_{1z} = q_1\}} g_z \Psi_{q_{2z}, q_{3z}, \chan_z}^{k_1}
\end{equation}
While this strategy increases the intermediate memory footprint by requiring the allocation of the $\Psi$ tensor, it removes the nested radial loops from the sparse execution path. This fused streaming approach is structurally advantageous for GPU architectures, where massive memory bandwidth can be leveraged to compute $\Psi$ in parallel, or when vectorizing the operator simultaneously over a large spatial grid.

\subsubsection{Angular-First Contraction (Pairing $\Gtens_{\mathbb{R}}$ first)}
To maximize the utilization of the CPU's high-speed L1 cache, the angular geometry can be aggregated before the dense physical tensor $\Rtens$ is applied. For a fixed target state $q_1$ and interaction channel $\chan$, we sweep through the sparse COO list to build a compact, intermediate angular tensor $\Phi$,
\begin{equation}
    \Phi_{k_2, k_3}^{(\chan, q_1)} = \sum_{(q_2, q_3) \in \Gtens_{\mathbb{R}}^{(\chan, q_1)}} g_z f_{k_2, q_2} f_{k_3, q_3}
\end{equation}
Because the radial basis size is typically small, e.g., $\kmax = 4$, the intermediate tensor $\Phi$ is simply a $5 \times 5$ matrix that fits within the CPU's local registers. Once constructed, this matrix is contracted against the dense physical tensor $\Rtens$,
\begin{equation}
    Q_{k_1, q_1} = \sum_{\chan=1}^{\numT} \sum_{k_2, k_3} \Rtens_{k_1, k_2, k_3}^{(\chan)} \Phi_{k_2, k_3}^{(\chan, q_1)}
\end{equation}
This strategy is optimal for local CPU execution. It decouples the sparse memory lookups from the heavy radial matrix-vector multiplications.

\subsection{Storage Complexity Analysis}
\label{sec:complexity}

Before empirically benchmarking the algorithms, it is instructive to derive the theoretical asymptotic scaling of the factorization in terms of both memory footprint and floating-point operations (FLOPS). By establishing the exact bounds of the discrete sums, the theoretical efficiency of the Wigner-Eckart decomposition can be quantified.

First, we define the global degrees of freedom for the spectral basis. Let $\kmax$ denote the maximum radial degree and $\lmax$ denote the maximum angular degree. The number of radial basis functions is $\numK = \kmax + 1$, and the number of angular basis functions is $\numQ = (\lmax + 1)^2$. The total number of spectral degrees of freedom per distribution function is therefore,
\begin{equation}
    \dof = \numK \numQ = (\kmax + 1)(\lmax + 1)^2.
\end{equation}

\begin{enumerate}
    \item \textbf{Storage Complexity of the Dense Cartesian Tensor:}
	The standard, un-factorized collision operator evaluates the full 3D Cartesian tensor $\Ctens_{\state_1, \state_2, \state_3}$. Because each of the three indices spans the entire basis $\dof$, the number of scalar elements in the dense tensor is,
\begin{equation}
    E_{\text{dense}} = \dof^3 = (\kmax + 1)^3 (\lmax + 1)^6.
\end{equation}
As $\lmax \to \infty$, the memory requirement of the dense tensor scales asymptotically as $\mathcal{O}(\kmax^3 \lmax^6)$.

 \item \textbf{Storage Complexity of the Factorized Operator:}
	The Wigner-Eckart factorization reduces this memory requirement by splitting the operator into the dense radial tensor $\Rtens$ and the sparse geometric routing table $\Gtens_{\mathbb{R}}$.

The dense physical tensor $\Rtens_{k_1, k_2, k_3}^{(\chan)}$ depends on the three radial indices and a single physical channel index $\chan$. This channel index enumerates all valid polar triplets $(l_1, l_2, l_3)$ that survive the geometric selection rules. The exact number of valid physical channels, $\numT$, is the sum of all triplets satisfying the triangle inequality and parity constraint,
\begin{equation}
    \numT = \sum_{l_1=0}^{\lmax} \sum_{l_2=0}^{\lmax} \sum_{\substack{l_3 = |l_1 - l_2| \\ (l_1+l_2+l_3) \text{ even}}}^{\min(l_1+l_2, \lmax)} 1 \approx \tfrac{1}{4} \lmax^3.
\end{equation}
The volume of this valid region scales proportionally to the volume of the angular truncation limit, yielding an asymptotic scaling of $\numT \sim \mathcal{O}(\lmax^3)$. The total number of elements in the dense physical tensor is $E_{\Rtens} = (\kmax + 1)^3 \numT$, which bounds its storage to $\mathcal{O}(\kmax^3 \lmax^3)$.

The geometric Gaunt tensor $\Gtens_{\mathbb{R}}$ contains the non-zero geometric transitions $\numG$ and is stored in a Coordinate (COO) format. This requires summing over all valid polar channels $\numT$, and then extracting the azimuthal modes $(m_1, m_2, m_3)$ that satisfy the zero-sum phase relation $m_1 = m_2 \pm m_3$,
\begin{equation}
    \numG = \sum_{(l_1, l_2, l_3) \in \numT} \sum_{m_1=-l_1}^{l_1} \sum_{m_2=-l_2}^{l_2} \sum_{m_3=-l_3}^{l_3} \delta_{m_1, m_2 \pm m_3}.
\end{equation}
Because the mode $m_3$ is constrained by the Kronecker delta, one degree of freedom is eliminated from the azimuthal summation, yielding an $\mathcal{O}(\lmax^2)$ multiplier. Combined with the $\mathcal{O}(\lmax^3)$ polar triplets, the number of non-zero geometric transitions scales as $\numG \sim \mathcal{O}(\lmax^5)$.

Therefore, the total storage complexity of the exact factorized representation is bounded by,
\begin{equation}
    \text{Storage}_{\text{Fact}} \sim \mathcal{O}(\kmax^3 \lmax^3) + \mathcal{O}(\lmax^5).
\end{equation}
\end{enumerate}

\subsection{Algorithmic Complexity Analysis}
While the asymptotic limits established in Section \ref{sec:complexity} dictate the global scaling of the factorized operator ($\mathcal{O}(\lmax^5)$), the physical wall-clock execution time is governed by the exact operation count (FLOPS) and the hardware-level memory access patterns. To elucidate the performance disparities observed in the benchmarks, we evaluate the algorithmic complexity of the four tensor contraction strategies. 

Let $\numG \sim \mathcal{O}(\lmax^5)$ denote the total number of non-zero geometric transitions. Furthermore, let $\numS$ denote the number of unique $(\chan, q_1)$ slices—the distinct combinations of physical interaction channels and target angular states. By definition, $\numS \le \numG$.

\begin{enumerate}
    \item \textbf{Dense Cartesian Contraction (Baseline):}
    The monolithic 3D tensor contraction evaluates the sum over the full spectral basis $\dof = \numK \numQ$. The operation count is bounded by,
    \begin{equation}
        \text{FLOPS}_{\text{Dense}} = \mathcal{O}(\kmax^3 \lmax^6).
    \end{equation}
    While iterating through a contiguous multi-dimensional array allows the CPU to prefetch memory into the L1 cache, the volume of zero-valued arithmetic operations dominates the execution time, establishing the upper computational boundary.

    \item \textbf{Standard Sparse Contraction:}
    This algorithm evaluates the discrete Wigner-Eckart sum as written, computing the full radial energy transfer for every non-zero geometric transition. 
    \begin{equation}
        \text{FLOPS}_{\text{Sparse}} = \mathcal{O}(\numG \kmax^3) \sim \mathcal{O}(\kmax^3 \lmax^5).
    \end{equation}
    Despite the mathematical reduction to $\mathcal{O}(\lmax^5)$, this implementation suffers from memory latency. As the algorithm increments through the 1D sparse array $z$, the angular indices $q_{2z}$ and $q_{3z}$ jump pseudo-randomly. This forces the CPU to constantly fetch non-contiguous segments of the distribution function $f$ from main RAM, stalling the arithmetic logic unit (ALU) due to L1 cache misses.

    \item \textbf{Radial-First Contraction:}
    By commuting the sum to evaluate the physical energy exchange prior to the geometric routing, this algorithm requires computing a dense intermediate physical state $\Psi_{z}^{k_1}$.
    \begin{equation}
        \text{FLOPS}_{\text{Radial}} = \mathcal{O}(\numG \kmax^3) + \mathcal{O}(\numG \kmax).
    \end{equation}
    While the asymptotic complexity matches the Standard Sparse approach, the memory profile is distinct. The algorithm must dynamically write and subsequently read the intermediate tensor $\Psi$, which contains $\numG \times \kmax$ elements. Because this tensor far exceeds the capacity of standard L1/L2 CPU caches, the execution becomes bound by the main memory bandwidth of the workstation.

    \item \textbf{Angular-First (Cache-Optimized Sliced):}
    The optimal contraction minimizes both memory latency and absolute floating-point operations. By sorting the data into $\numS$ unique slices, the algorithm first builds the intermediate angular matrix $\Phi$ (requiring $\kmax^2$ operations per transition), and subsequently contracts $\Phi$ against the dense physical tensor $\Rtens$ (requiring $\kmax^3$ operations per slice). The total operation count is,
    \begin{equation}
        \text{FLOPS}_{\text{Angular}} = \mathcal{O}(\numG \kmax^2) + \mathcal{O}(\numS \kmax^3).
    \end{equation}
    Because the number of unique slices is less than the number of total transitions ($\numS \le \numG$), factoring the $\kmax^3$ radial loop out of the sparse geometric loop reduces the total number of required arithmetic operations compared to the Standard Sparse and Radial-First algorithms. 
    
    Coupled with the fact that the $\kmax \times \kmax$ matrix $\Phi$ resides within the CPU's local registers, this algorithm achieves theoretical minimum arithmetic complexity while executing contiguous, cache-friendly dense matrix multiplications. This dual optimization is the origin of the acceleration observed in the empirical benchmarks.
\end{enumerate}

%% file: section5.tex
\section{Results}
This section evaluates the factorized collision operator in three parts. First, we demonstrate the spectral convergence of the singular quadrature scheme for both smooth and non-analytic kernels. Second, we validate the physical fidelity of the assembled tensors: for Maxwell molecules, we confirm the exact conservation of macroscopic invariants, adherence to the discrete H-theorem, and recovery of the classical Wang Chang-Uhlenbeck spectrum; for hard spheres, we verify the capture of nonlinear energy cascades and the recovery of the infinite-order Chapman-Enskog transport coefficients. Finally, we benchmark memory footprint and execution time of the factorized contractions against the Cartesian baseline.

\subsection{Spectral Convergence of the Singular Quadrature}
To quantify the accuracy of the singular quadrature scheme developed in Section \ref{sec:ingularquad}, we evaluate the numerical integration error for both Hard Sphere ($\gamma=1$) and Maxwell molecule ($\gamma=0$) gases. Because the non-analytic collision kernel ($u^\gamma$) and the mapped Gaussian envelope prevent standard algebraic exactness, we measure the convergence against an over-resolved numerical reference tensor $\Rtens_{\text{ref}}$, with a quadrature padding of $N_{\text{pad}} = 64$. 

For a fixed spectral basis ($K_{\max}=4, L_{\max}=4$), the tensor is assembled using padding values of $N_{\text{pad}} \in \{0, 2, 4, 8, 16, 32\}$. This padding is applied strictly to the coupled kinematic coordinates associated with the non-polynomial terms. Because the rotational scattering angles ($\chi, \epsilon$) govern only finite polynomial expansions via the spherical harmonics, their theoretical exactness bounds remain sufficient. The global tensor error is computed via the relative $\ell_\infty$ maximum norm:
\begin{equation}
    e = \frac{||\Rtens_N - \Rtens_{\text{ref}}||_{\ell_\infty}}{||\Rtens_{\text{ref}}||_{\ell_\infty}}.
\end{equation}

\begin{figure}[htbp]
    \centering
    \includegraphics[width=0.6\textwidth]{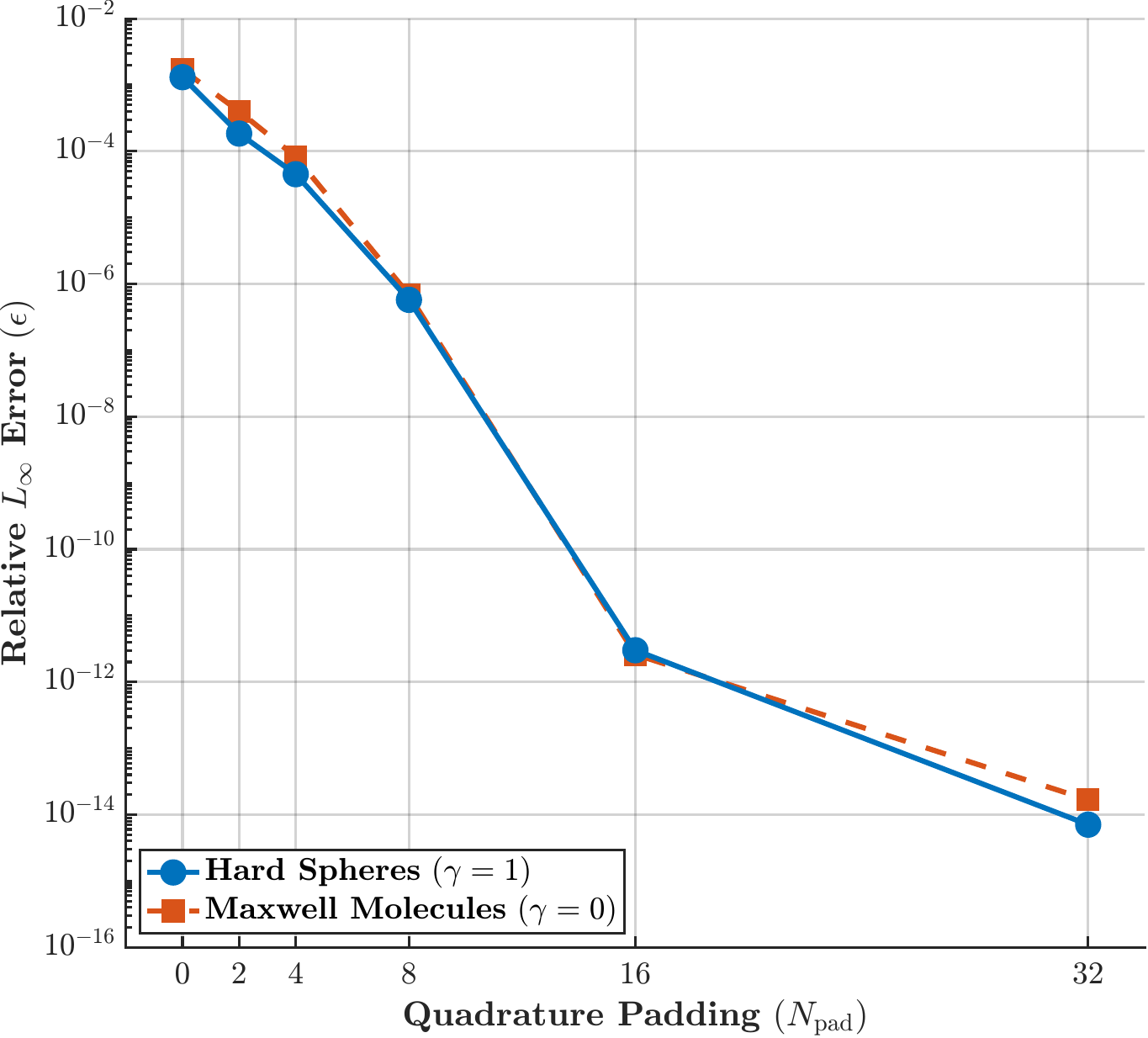}
    \caption{Relative $\ell_\infty$ tensor error for Hard Sphere ($\gamma=1$) and Maxwell molecule ($\gamma=0$) gases ($K_{\max}=4, L_{\max}=4$) as a function of the 2D Duffy quadrature padding. Both interaction potentials exhibit identical exponential decay to the double-precision machine epsilon.}
    \label{fig:quadrature_convergence}
\end{figure}

As demonstrated in Figure \ref{fig:quadrature_convergence}, both interaction limits exhibit identical exponential convergence. The overlapping decay rates demonstrate that the 2D Duffy transformation regularizes the cone singularity of the Hard Sphere kernel, allowing the Gauss-Legendre tensor-product grid to evaluate the singular integrals with the same spectral accuracy as the smooth Maxwell molecule kernel.

\subsection{Maxwell Molecules: Analytical Benchmarks}
We verify the core implementation by testing the operator against analytical results, including the transient Bobylev-Krook-Wu (BKW) solution and the Wang Chang-Uhlenbeck eigenvalue spectrum. We also confirm adherence to the discrete H-theorem and the conservation of macroscopic invariants.

\subsubsection{The Bobylev-Krook-Wu (BKW) Solution}
To verify the accuracy of the numerical framework, we evaluate the nonlinear collision operator against the Bobylev-Krook-Wu (BKW) solution~(\cite{bobylev:1975}). The BKW solution is a rare analytical benchmark for the Boltzmann equation, describing the isotropic relaxation of a non-equilibrium gas of Maxwellian molecules.

\begin{figure}[htbp]
    \centering
    \begin{subfigure}{\linewidth}
        \centering
        \includegraphics[width=0.75\linewidth]{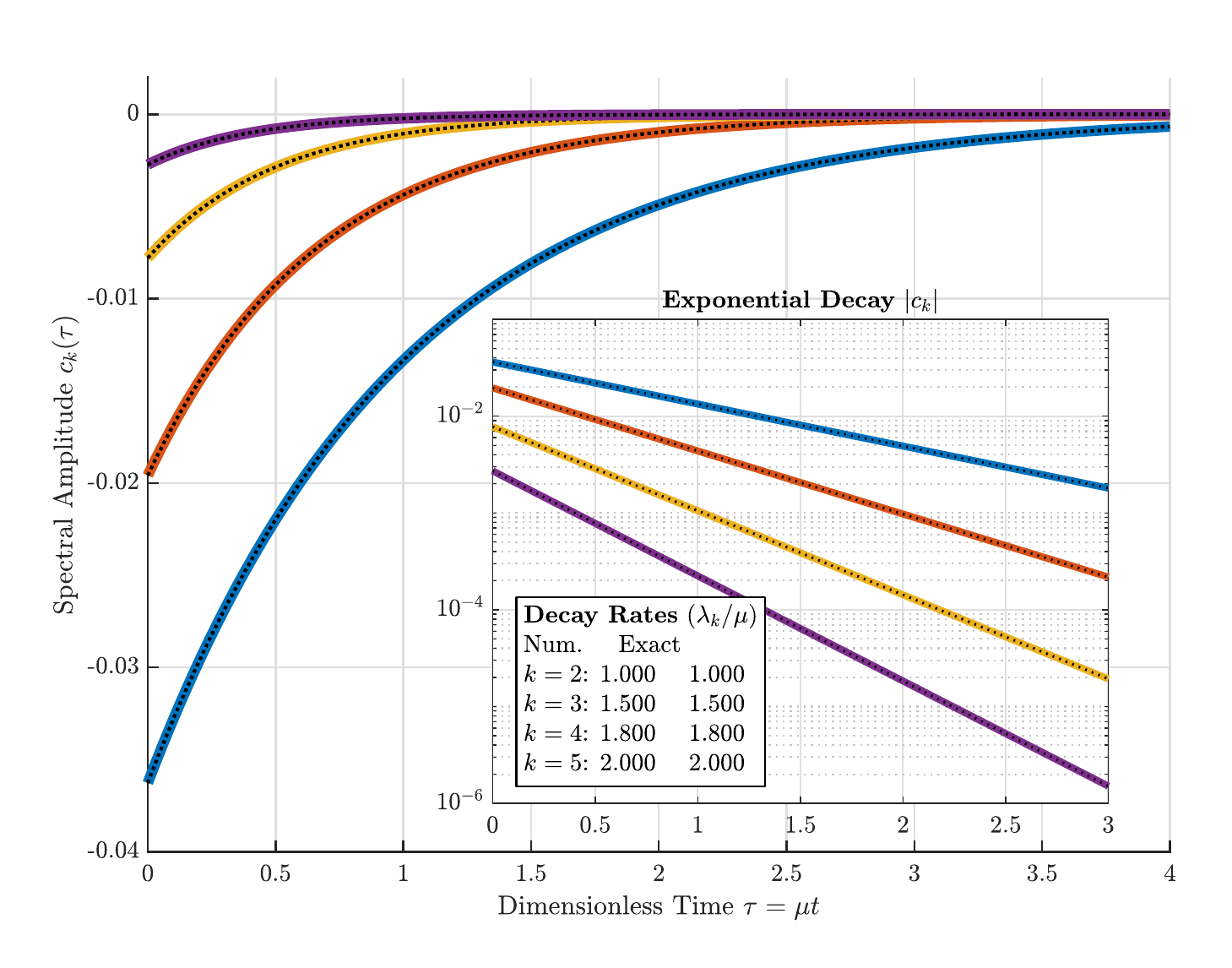}
        \caption{Nonlinear spectral relaxation and asymptotic decay rates.}
        \label{fig:bkw_relax}
    \end{subfigure}
    \vspace{1.5em}
    \begin{subfigure}{\linewidth}
        \centering
        \includegraphics[width=0.75\linewidth]{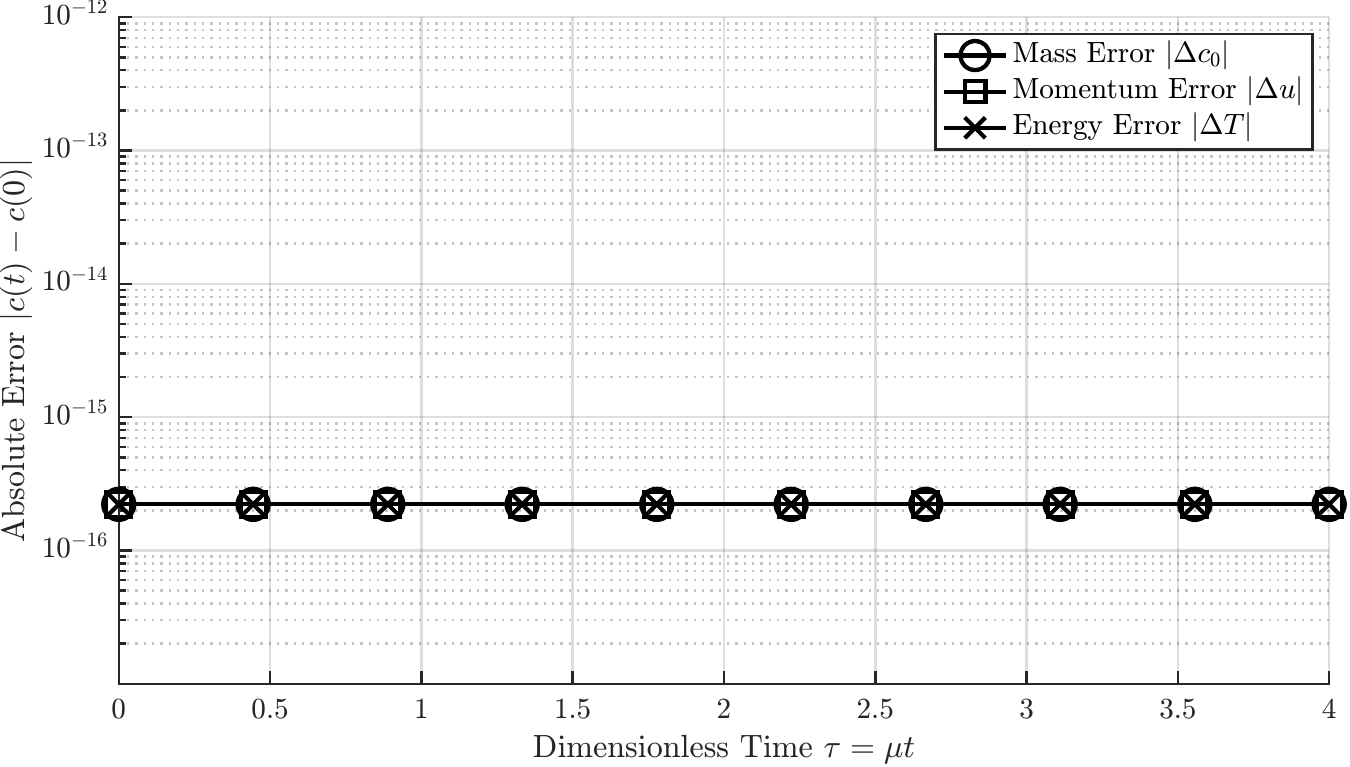}
        \caption{Machine-precision conservation of collision invariants.}
        \label{fig:bkw_cons}
    \end{subfigure}
    \caption{Validation of the nonlinear collision operator against the Bobylev-Krook-Wu (BKW) solution. (a) Time evolution of the isotropic spectral amplitudes $c_k(\tau)$ for $k \ge 2$. The inset confirms that the numerical relaxation captures the theoretical Wang Chang-Uhlenbeck (WCU) eigenvalues in the linear regime. (b) Absolute error in the macroscopic mass and energy invariants, demonstrating numerical conservation down to machine precision.}
    \label{fig:bkw_validation}
\end{figure}

Figure \ref{fig:bkw_validation}(a) displays the time evolution of the isotropic spectral amplitudes during this relaxation process. As the gas transitions from a nonlinear state into the linear regime, the extracted decay rates of the individual modes match the theoretical \cite{wangchang:1951} (WCU) eigenvalues. Furthermore, Figure \ref{fig:bkw_validation}(b) confirms that this dynamic energy exchange does not introduce numerical artifacts; the macroscopic mass and energy invariants are conserved to machine precision.

\subsubsection{Eigenvalue Spectrum and the Discrete H-Theorem}

Beyond simulating transient relaxation, a spectral operator must possess stable structural properties, including adherence to the discrete H-theorem and the exact preservation of collision invariants.

\begin{figure}[htbp]
    \centering
    \includegraphics[width=1\textwidth]{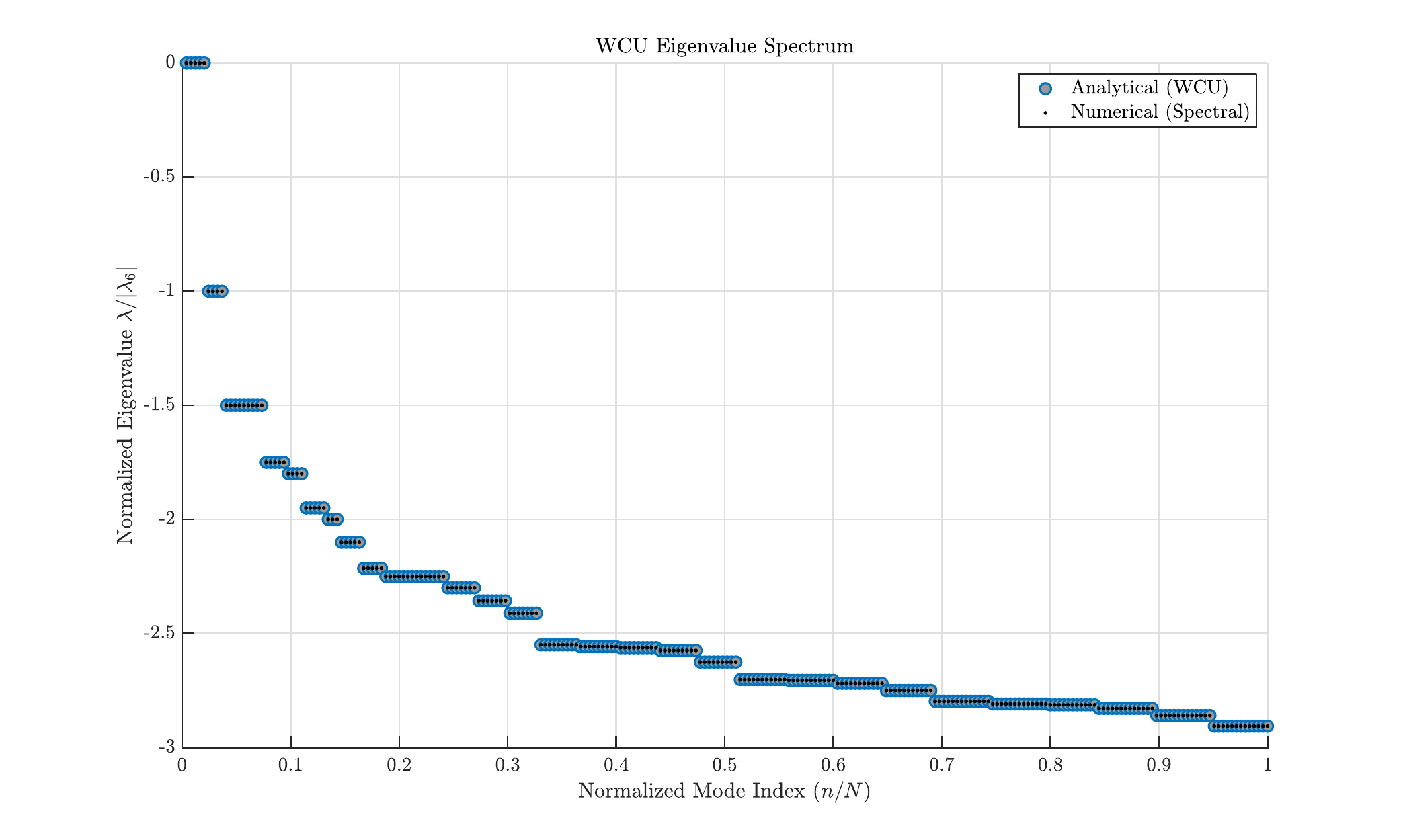}
    \caption{Wang Chang-Uhlenbeck (WCU) eigenvalue spectrum of the linearized collision operator for Maxwell molecules ($K_{\max}=4, L_{\max}=6$). The numerical eigenvalues (black dots) match the analytical predictions (grey circles). The five macroscopic collision invariants are isolated at $\lambda = 0$, while the relaxation modes display a characteristic "staircase" structure due to mode degeneracy, clustering near the lower bound of \num{-3}.}
    \label{fig:wcu_spectrum}
\end{figure}

We analyze the eigenvalue spectrum of the nonlinear collision operator linearized around a Maxwellian equilibrium. For Maxwell molecules, the exact eigenvalues are known analytically through Wang Chang-Uhlenbeck (WCU) theory \cite{wangchang:1951}. As shown in Figure \ref{fig:wcu_spectrum} and detailed in Table \ref{tab:wcu_data}, the numerical eigenvalues generated by the factorized operator match the analytical WCU predictions. The numerical spectrum is strictly non-positive, satisfying the discrete H-theorem. The five macroscopic collision invariants (mass, momentum, and energy) are captured exactly at $\lambda = 0$, while all other modes drive the system toward equilibrium with negative decay rates. These relaxation modes capture the exact symmetric mode degeneracies, resulting in a characteristic "staircase" structure that clusters near the theoretical lower bound. Furthermore, the maximum imaginary component of the numerical eigenvalues is bounded by machine precision, confirming that the spectral approximation rigorously preserves the self-adjoint symmetry of the continuous collision operator.

\begin{table}[htbp]
    \centering
    \caption{Numerical verification of the linearized collision operator against WCU theory. Eigenvalues are normalized by the first physical relaxation mode $|\lambda_6|$. The grouping highlights the capture of symmetric degeneracies.}
    \label{tab:wcu_data}
    \begin{tabular}{@{}lccccl@{}}
        \toprule
        Mode Index $n$ & Num. Ratio & Ana. Ratio & Abs. Diff. & Interpretation \\ \midrule
        1 -- 5   &  0.0000 &  0.0000 & \num{0.00e00} & Collision Invariants \\
        6 -- 9   & -1.0000 & -1.0000 & \num{< 3e-15} & Fundamental Relaxation \\
        10 -- 18 & -1.5000 & -1.5000 & \num{< 2e-14} & Higher-order Coupling \\
        19 -- 23 & -1.7500 & -1.7500 & \num{< 2e-14} & Higher-order Coupling \\
        24 -- 27 & -1.8000 & -1.8000 & \num{< 1e-14} & Higher-order Coupling \\ 
        \bottomrule
    \end{tabular}
\end{table}

\subsubsection{Galilean Invariance}
Finally, Table \ref{tab:galilean_invariance} demonstrates the Galilean invariance of the discrete operator. When the equilibrium Maxwellian is shifted by a bulk velocity $\vmean$, the discrete basis experiences an increase in spectral truncation error. However, because the conservation properties are embedded in the geometric null space of the operator, the macroscopic mass, momentum, and energy remain conserved to machine precision regardless of the truncation error.

\begin{table}[htbp]
    \centering
    \caption{Galilean invariance test for a shifted Maxwellian distribution ($K_{\max}=4, L_{\max}=6$). As the bulk velocity $\vmean$ increases, the distribution shifts away from the origin, leading to an increase in the global spectral truncation error ($L_2$ norm). However, despite the spectral truncation at higher velocities, the macroscopic collision invariants (mass, momentum, and energy) remain conserved to machine precision.}
    \label{tab:galilean_invariance}
    \begin{tabular}{@{}lccc@{}}
        \toprule
        Bulk Velocity $\vmean$ & Truncation Error ($L_2$) & Conservation Error & Configuration \\ \midrule
        $[0.0, 0, 0]$ & \num{4.77e-15} & \num{0.00e+00} & Stationary (Exact) \\
        $[0.1, 0, 0]$ & \num{5.14e-11} & \num{0.00e+00} & Shifted (Truncated) \\
        $[0.3, 0, 0]$ & \num{3.39e-07} & \num{0.00e+00} & Shifted (Truncated) \\
        $[0.6, 0, 0]$ & \num{8.88e-05} & \num{0.00e+00} & Shifted (Truncated) \\ \bottomrule
    \end{tabular}
\end{table}

\subsection{Hard Spheres: Asymptotic Transport Limits}
\label{sec:results_chapman_enskog}
We extend the validation to velocity-dependent collision models. By simulating an anisotropic Hard Sphere gas, we demonstrate that the numerical operator captures nonlinear energy cascades and recovers the theoretical infinite-order Chapman-Enskog viscosity limits.

\subsubsection{Infinite-Order Chapman-Enskog Inversion}
While the BKW solution validates the operator for Maxwellian molecules, simulating realistic gas flows requires the operator to handle velocity-dependent collision kernels, such as the Hard Sphere model. Because no analytical solution exists for the distribution function of a Hard Sphere gas relaxing to equilibrium, validation relies on recovering the macroscopic transport coefficients defined by the Chapman-Enskog expansion (\cite{chapman:1970, ferziger:1972}).

\paragraph{\textbf{Continuous Governing Equation.}} In classical Chapman-Enskog theory, the shear viscosity is derived from the first-order perturbation $f^{(1)}$ to the local Maxwellian equilibrium $f^{(0)}$. This perturbation is governed by the linearized collision operator $\mathcal{L}$ driven by the traceless symmetric strain rate tensor:
\begin{equation}
    f^{(0)} \left( \vpre \otimes \vpre - \frac{1}{3}\vmag^2 \mathbf{I} \right) : \nabla \vmean = \mathcal{L}(f^{(1)}),
\end{equation}
where $\vmean$ is the macroscopic fluid velocity, $\mathbf{I}$ is the identity tensor, and $\vpre$ denotes the microscopic particle velocity with magnitude $\vmag = |\vpre|$.

Traditionally, calculating the viscosity requires expanding $f^{(1)}$ into an infinite series of Sonine polynomials and computing a dense matrix of continuous bracket integrals \citep{pekeris:1957}, which becomes analytically prohibitive beyond the first few terms. The Wigner-Eckart factorization bypasses this approach.

\paragraph{\textbf{Spectral Block-Diagonalization.}} Because the macroscopic source term on the left-hand side is a traceless symmetric rank-2 tensor, it maps exclusively to the spherical harmonic degree $L=2$. Consequently, when projected onto the discrete Wigner-Eckart spectral basis, the source translates to a sparse vector $\mathbf{s}$ containing non-zero entries exclusively for modes where $l=2$.

To solve for the discrete perturbation coefficients $\delta \mathbf{c}$, we construct the spectral matrix representation of the linearized collision operator. Evaluated at the absolute equilibrium state $\mathbf{c}_{\text{eq}}$, where only the isotropic ground mode $c_{0,0,0} = 1$ is non-zero, this discrete linear operator $\mathbf{L}$ is given by,
\begin{equation}
    \mathbf{L}_{\state_1, \statepert} = \frac{\partial Q_{\state_1}}{\partial c_{\statepert}} \Bigg|_{\mathbf{c}_{eq}} = \Ctens_{\state_1, (0,0,0), \statepert} + \Ctens_{\state_1, \statepert, (0,0,0)}
\end{equation}

Substituting the Wigner-Eckart factorization $\Ctens_{\state_1, \state_2, \state_3} = \Rtens_{k_1 k_2 k_3}^{l_1 l_2 l_3} \; \Gtens_{m_1 m_2 m_3}^{l_1 l_2 l_3}$ into the Jacobian block-diagonalizes the system. Because the background Maxwellian is isotropic, it forces the intermediate angular degrees to be $l_2 = 0$ and $m_2 = 0$. Applying the Clebsch-Gordan selection rules enforced by the Gaunt tensor $\Gtens$, the non-zero transitions are bounded by:
\begin{enumerate}
    \item \textbf{Triangle Inequality:} $|0 - l_{\statepert}| \le l_1 \le 0 + l_{\statepert} \implies l_1 = l_{\statepert}$
    \item \textbf{Azimuthal Balance:} $m_1 = 0 + m_{\statepert} \implies m_1 = m_{\statepert}$
\end{enumerate}

This geometric constraint dictates that an isotropic background cannot alter the angular momentum of a colliding perturbation. The Jacobian matrix evaluated at equilibrium block-diagonalizes with respect to both the polar degree $l$ and azimuthal mode $m$. The continuous integro-differential equation reduces to a finite, decoupled linear system governed exclusively by the $\kmax \times \kmax$ radial block for $L=2$:
\begin{equation}
    \mathbf{L}^{(L=2)} \cdot \delta \mathbf{c}_{L=2} = \mathbf{s}_{L=2}
\end{equation}

By directly inverting this $\kmax \times \kmax$ block, the spectral operator yields the shear viscosity to arbitrary order without evaluating continuous bracket integrals. The $k=0$ component of the solution vector corresponds to the first-order Chapman-Enskog approximation, while the higher-order radial modes $k \ge 1$ incorporate the higher-order Sonine polynomial corrections.

\paragraph{\textbf{Numerical Validation.}} To validate this numerical implementation, the factorized operator was evaluated for a Hard Sphere gas ($\gamma = 1$) and the corresponding $\mathbf{J}^{(L=2)}$ block was inverted. Table \ref{tab:viscosity_correction} compares the numerical viscosity correction factors $f_\mu^{(K)} = \mu^{(K)} / \mu^{(1)}$ against the analytical limits derived from classical bracket integrals.

\begin{table}[h!]
\centering
\caption{Numerical verification of the infinite-order Chapman-Enskog shear viscosity correction factors $f_\mu^{(K)}$ for a Hard Sphere gas. The spectral Wigner-Eckart inversion asymptotes to the infinite-order theoretical limit without evaluating continuous bracket integrals.}
\label{tab:viscosity_correction}
\begin{tabular}{@{}lccc@{}}
\toprule
\textbf{Truncation Limit} & \textbf{Matrix Size} & \textbf{Numerical $f_\mu$} & \textbf{Exact Analytical $f_\mu$} \\ \midrule
$\kmax = 0$             & $1 \times 1$         & 1.000000                   & 1.000000 (1st-Order)              \\
$\kmax = 1$             & $2 \times 2$         & 1.014851                   & 1.014851 (2nd-Order)              \\
$\kmax = 2$             & $3 \times 3$         & 1.015879                   & --                                \\
$\kmax = 3$             & $4 \times 4$         & 1.016006                   & --                                \\
$\kmax = 4$             & $5 \times 5$         & 1.016028                   & --                                \\ \midrule
\textbf{Theory ($\infty$)}& $\infty \times \infty$& --                         & \textbf{1.016034 (Inf-Order)}     \\ \bottomrule
\end{tabular}
\end{table}

For $\kmax = 1$, the spectral inversion recovers the second-order analytical fraction $205/202 \approx 1.014851$. As $\kmax$ increases, the numerical viscosity asymptotes to the theoretical infinite-order ceiling of $1.016034$. Both of these analytical limits match the classical Sonine polynomial sequence computed by \cite{pekeris:1957}. This confirms that the Wigner-Eckart factorization evaluates the velocity-dependent physics of the Hard Sphere kernel down to machine precision, capturing the higher-order cross-couplings of the continuous Boltzmann integral without numerical artifacts.

\subsubsection{Transient Anisotropic Stress Relaxation}
\begin{figure}[h!]
    \centering
    \begin{subfigure}{0.75\textwidth}
        \centering
        \includegraphics[width=\linewidth]{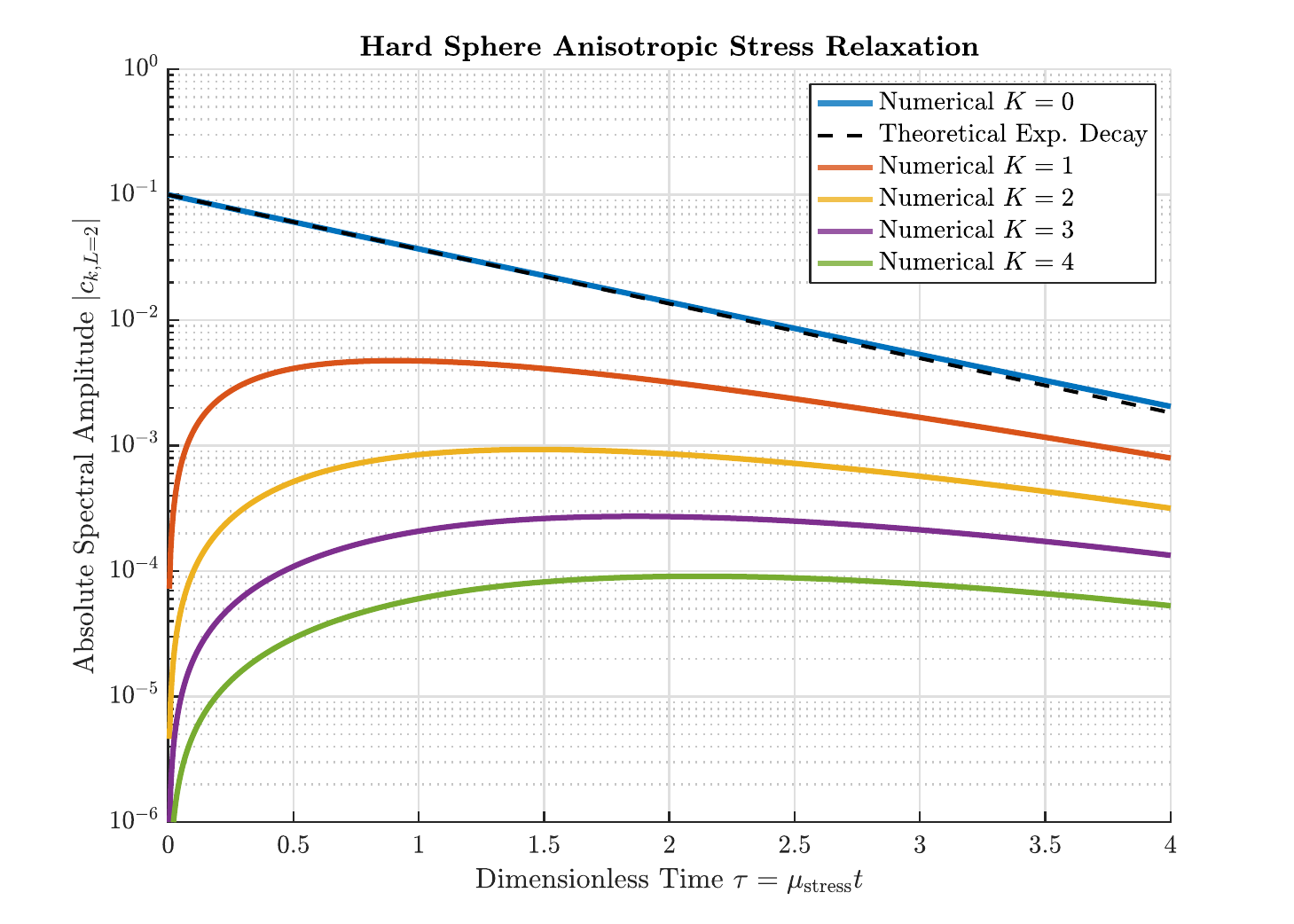}
        \caption{Anisotropic Stress Relaxation}
        \label{fig:stress_relax}
    \end{subfigure}\hfill
    \begin{subfigure}{0.75\textwidth}
        \centering
        \includegraphics[width=\linewidth]{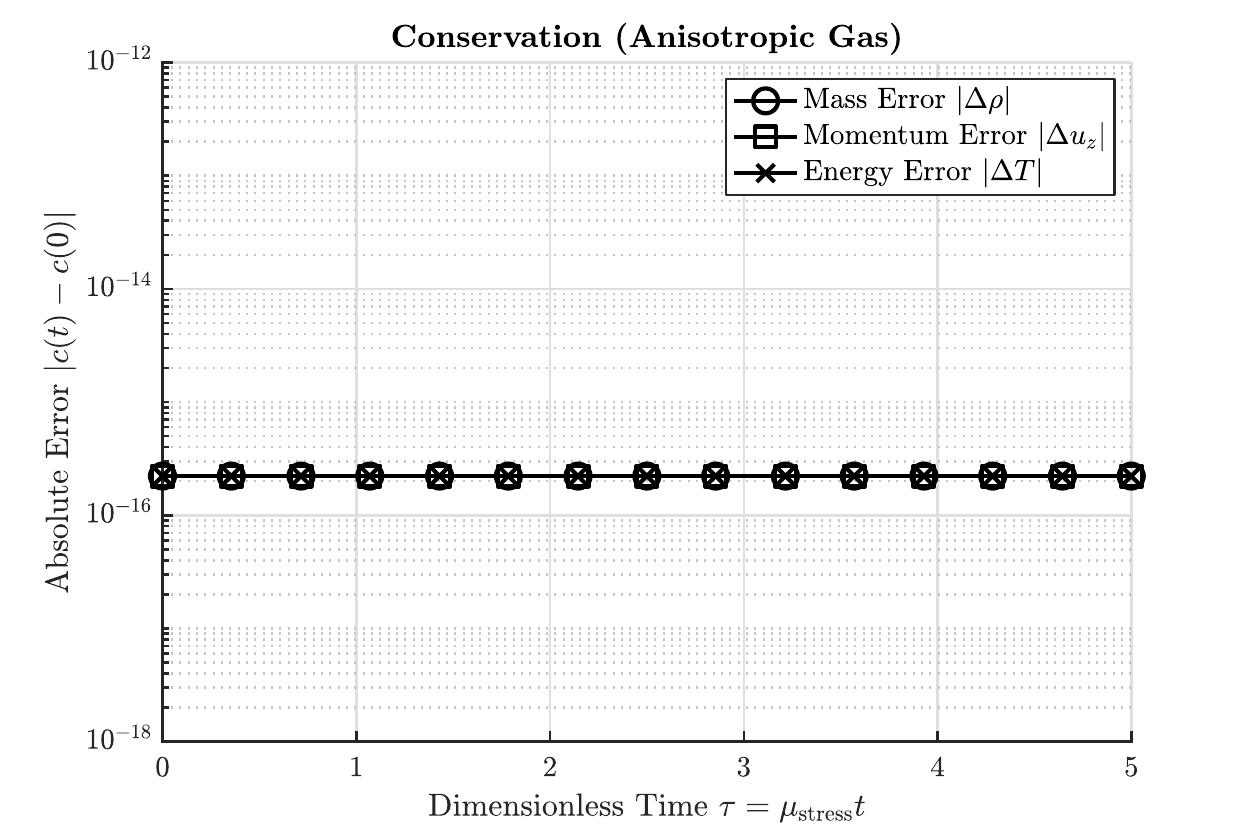}
        \caption{Conservation of Invariants}
        \label{fig:stress_conserv}
    \end{subfigure}
    \caption{Transient relaxation of a non-equilibrium macroscopic shear stress in a Hard Sphere gas ($\gamma = 1$). (a) The primary $K=0$ shear mode tracks the theoretical Chapman-Enskog exponential decay dictated by the $f_\mu$ viscosity correction. The nonlinear collision operator scatters energy into higher-order radial polynomials ($K>0$) before thermalization. (b) The macroscopic collision invariants (mass, momentum, and energy) remain conserved to machine precision throughout the nonlinear anisotropic relaxation.}
    \label{fig:transient_stress}
\end{figure}

To confirm the temporal stability of these extracted transport coefficients, the factorized operator was integrated over continuous time using an explicit Runge-Kutta scheme. The gas was initialized as a Maxwellian superimposed with a macroscopic anisotropic stress perturbation in the primary $L=2, K=0$ mode. As shown in Figure \ref{fig:transient_stress}(a), the primary numerical shear mode tracks the theoretical exponential decay dictated by the Chapman-Enskog viscosity correction. Furthermore, the operator stably captures the nonlinear energy cascade into higher-order radial polynomials ($K>0$) while preserving the macroscopic collision invariants to machine precision (Figure \ref{fig:transient_stress}b).

\subsection{Computational Performance}
We quantify the hardware advantages of the factorization. By benchmarking both memory footprint and execution time, we demonstrate how the compressed geometric routing overcomes the bottlenecks of standard dense Cartesian methods.

\subsubsection{Execution Time}
The primary advantage of the Wigner-Eckart factorization is its ability to overcome the computational and memory bottlenecks that limit spectral Boltzmann solvers.

\begin{figure*}[h!]
    \centering
    \begin{subfigure}[b]{0.49\textwidth}
        \centering
        \includegraphics[width=\linewidth]{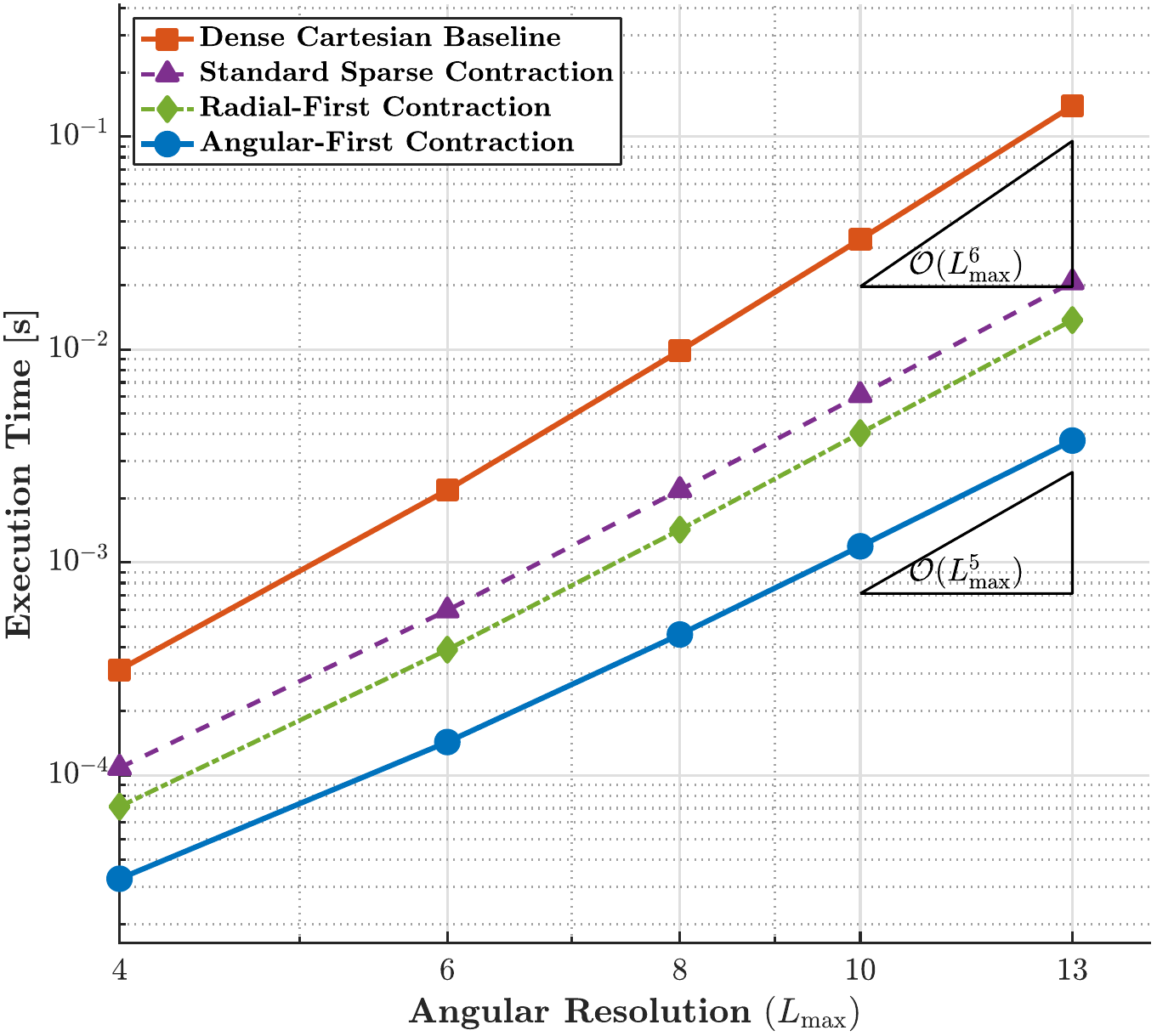}
        \caption{Absolute execution time [s]}
        \label{fig:exec_absolute}
    \end{subfigure}
    \hfill
    \begin{subfigure}[b]{0.49\textwidth}
        \centering
        \includegraphics[width=\linewidth]{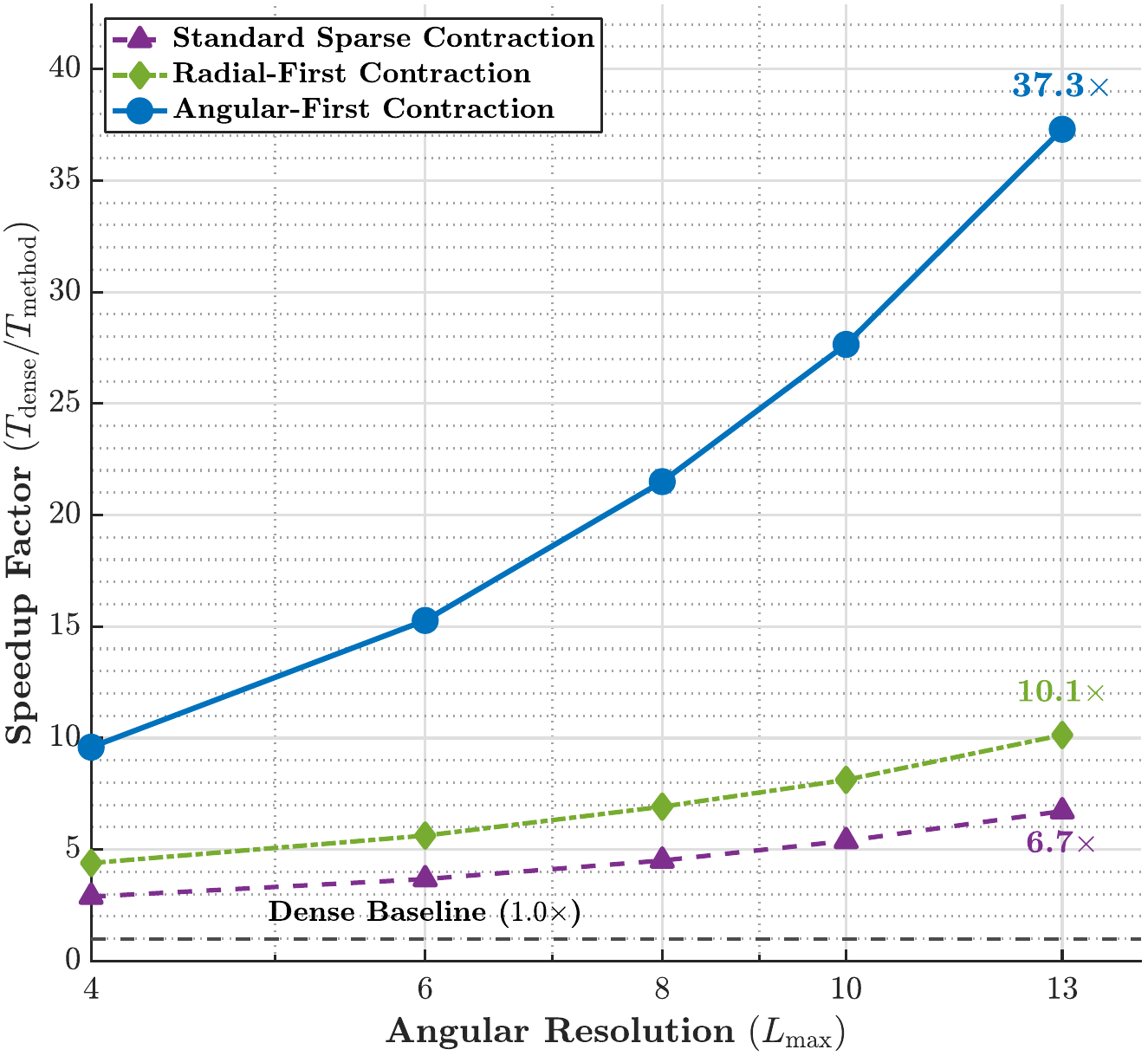}
        \caption{Relative speedup ($T_{\text{dense}}/T_{\text{method}}$)}
        \label{fig:exec_relative}
    \end{subfigure}
    \caption{Computational complexity and hardware efficiency analysis of the Boltzmann collision operator evaluation ($\kmax=4$). (a) Absolute execution time scaling demonstrating the theoretical $\mathcal{O}(\lmax^6)$ dense baseline versus the $\mathcal{O}(\lmax^5)$ factorized algorithms. (b) The linear speedup factor highlights how the cache-optimized sliced contraction minimizes memory latency, achieving a $37.2\times$ acceleration over the standard baseline approach at $\lmax=13$.}
    \label{fig:execution_results}
\end{figure*}

Figure \ref{fig:execution_results}(a) compares the absolute execution time of the factorized algorithms against a standard dense Cartesian baseline. As shown in Figure \ref{fig:execution_results}(b), the cache-optimized angular-first contraction minimizes memory latency, yielding a $37.2\times$ acceleration over the baseline at an angular resolution of $\lmax=13$.

\begin{figure*}[h!]
    \centering
    \begin{subfigure}[b]{0.49\textwidth}
        \centering
        \includegraphics[width=\linewidth]{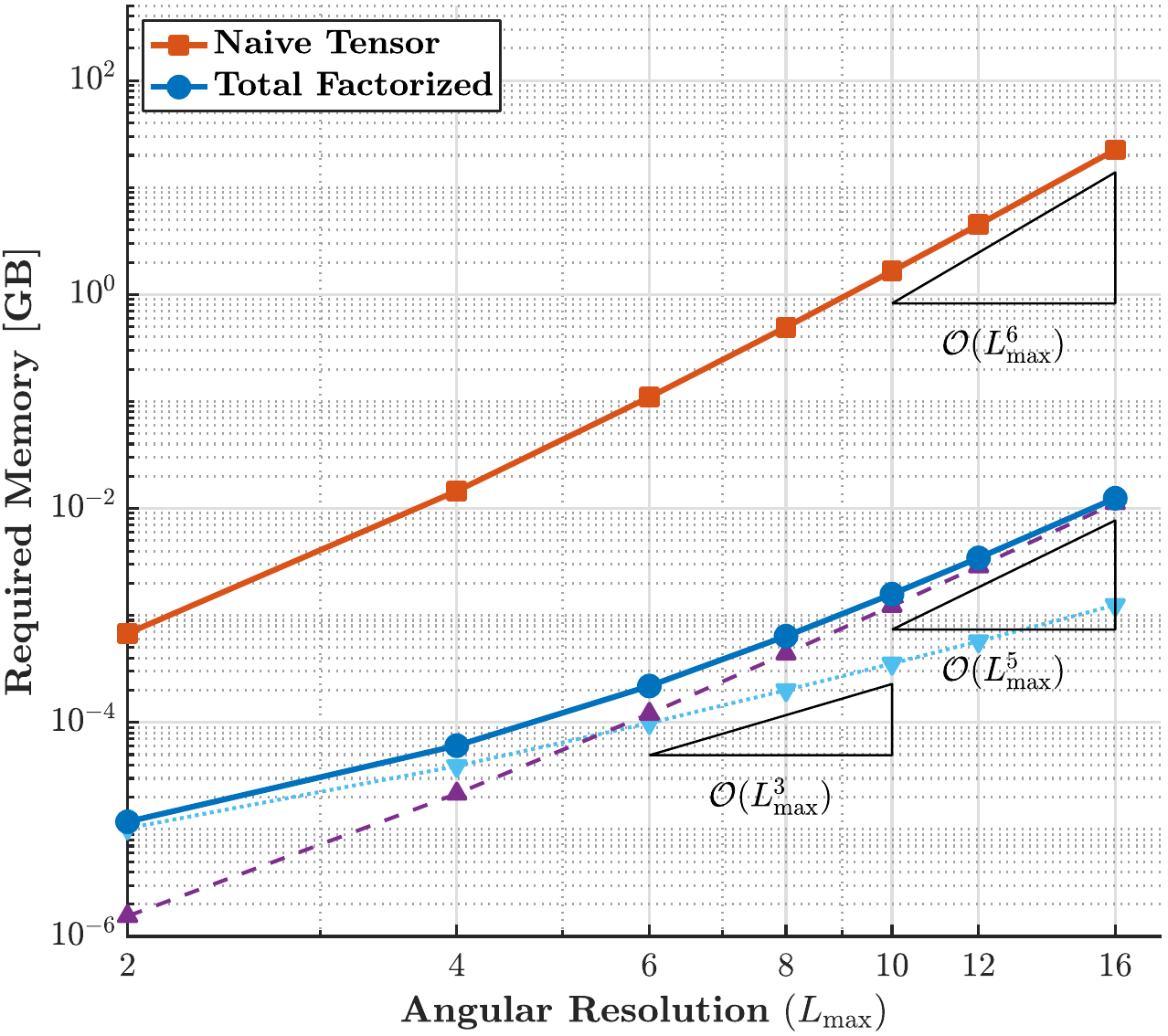}
        \caption{Absolute memory footprint [GB]}
        \label{fig:mem_absolute}
    \end{subfigure}
    \hfill
    \begin{subfigure}[b]{0.49\textwidth}
        \centering
        \includegraphics[width=\linewidth]{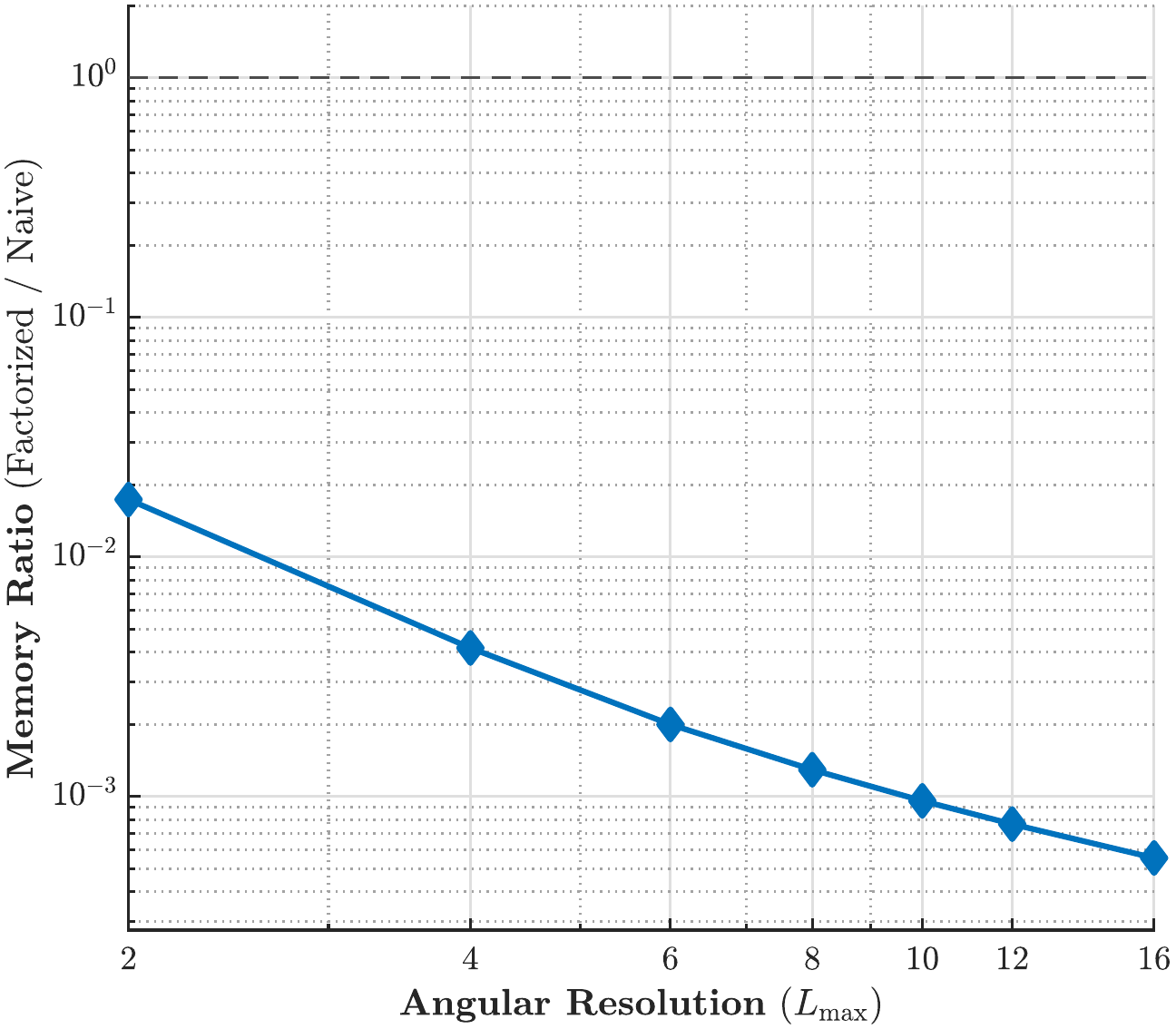}
        \caption{Relative efficiency (Factorized/Naive)}
        \label{fig:mem_relative}
    \end{subfigure}
    \caption{Memory complexity analysis of the Wigner-Eckart factorized collision operator ($\kmax=4$). 
    (a) Comparison of the memory footprint between the naive Cartesian tensor ($\mathcal{O}(\lmax^6)$) and the proposed decomposition ($\mathcal{O}(\lmax^5)$ and $\mathcal{O}(\lmax^3)$ regimes). 
    (b) Normalized memory ratio demonstrating a three-order-of-magnitude reduction in storage requirements for high angular resolutions.}
    \label{fig:memory_results}
\end{figure*}

\subsubsection{Memory Complexity and Storage Reduction}
This acceleration is enabled by a reduced memory footprint. As illustrated in Figure \ref{fig:memory_results} and Table \ref{tab:memory_scaling}, the memory required to store the dense Cartesian tensor scales asymptotically as $\mathcal{O}(\lmax^6)$. By storing the macroscopic geometry in a Coordinate (COO) format, the factorization reduces storage requirements by three orders of magnitude at high angular resolutions.

\begin{table}[tbp]
\centering
\caption{Memory complexity and storage efficiency for the Boltzmann collision tensor ($\kmax=4$). The factorized memory includes both the dense collision tensor $\Rtens$ and the sparse geometric transitions $\Gtens$ in COO format.}
\label{tab:memory_scaling}
\begin{tabular}{@{}l r r r r S[table-format=1.2e-1]@{}}
\toprule
$\lmax$ & DOFs & Gaunt Nonzeros & Naive [GB] & Fact. [GB] & {Ratio (Factorized / Naive)} \\
\midrule
6  & 245  & 6,460   & 0.11 & 0.22e-03 & 2.00e-03 \\
8  & 405  & 23,621  & 0.50 & 0.66e-03 & 1.29e-03 \\
10 & 605  & 65,913  & 1.65 & 1.62e-03 & 9.59e-04 \\
12 & 845  & 154,330 & 4.50 & 3.53e-03 & 7.67e-04 \\
16 & 1,445 & 601,569 & 22.50 & 12.75e-03 & 5.54e-04 \\
\bottomrule
\end{tabular}
\end{table}

%% file: section6.tex
\section{Conclusion}
We presented a conservative spectral method for the bilinear Boltzmann collision operator based on a geometric dimensional reduction. By rigidly rotating the laboratory frame to align with the colliding pair and integrating over the $\mathrm{SO}(3)$ rotation group, we decoupled the physical scattering dynamics from the macroscopic angular geometry. 

The implementation and performance of this numerical framework are defined by five methodological contributions:
\begin{itemize}
    \item \textbf{Geometric Factorization:} Reducing the 8-dimensional Cartesian collision integral to a 5-dimensional kinematic core, establishing the Wigner-Eckart factorization without abstract algebraic transformations.
    \item \textbf{Spectral Quadrature:} Resolving the non-analytic velocity-dependent collision kernel through a 2D Duffy transformation, enabling Gauss-type quadrature rules to achieve exponential convergence.
    \item \textbf{Storage Compression:} Exploiting Clebsch-Gordan selection rules to store the macroscopic geometry in a sparse Coordinate (COO) format, reducing the operator memory footprint by three orders of magnitude.
    \item \textbf{Invariant Conservation:} Embedding the physical null-space directly into the discrete tensor by zeroing specific entries, guaranteeing conservation of mass, momentum, and energy to machine precision.
    \item \textbf{Execution Acceleration:} Utilizing a cache-optimized, angular-first tensor contraction strategy that decouples sparse memory routing from dense arithmetic, yielding a 37-fold execution speedup over standard formulations.
\end{itemize}

The numerical framework was validated against the analytical Bobylev-Krook-Wu solution and the Wang Chang-Uhlenbeck eigenvalue spectrum for Maxwell molecules, as well as the infinite-order Chapman-Enskog viscosity limits for hard spheres. By eliminating the $\mathcal{O}(L_{\max}^6)$ memory and arithmetic bottlenecks inherent to dense Cartesian methods, this factorization provides a robust and efficient foundation for high-resolution three-dimensional deterministic kinetic simulations.

Interesting directions for future research include extending the methodology to first-principles polyatomic collision models, such as the Waldmann-Snider equation, and applying the framework to the spatially inhomogeneous setting.

%% file: appendix1.tex
\section{Cascaded Evaluation and Sum Factorization of the Kinematic Core} \label{app:1}

Evaluating the discrete five-dimensional kinematic core directly on the tensor-product grid introduces an $\mathcal{O}(N_{\cmenergy} N_{\fracrel} N_{\duffy} N_{\defangle} N_{\aziangle})$ arithmetic complexity per physical channel. To optimize the evaluation, the numerical integration is executed as a sequence of block-decoupled tensor contractions. By recursively applying Fubini's theorem, we cascade the integral from the inside out, leveraging the structure of the Duffy patches to sum-factorize the physical components.

\subsection{Innermost Scattering Manifold}

We arrange the discrete grid generation such that the scattering angles $\defangle$ and $\aziangle$ form the innermost integration loops. This ensures that the pre-collision relative speed $\urel$ remains constant during the scattering evaluation, allowing the pre-collision loss term to be subtracted from the azimuthally-integrated gain term before the collision kernel is applied. Because the loss term is uncoupled from the scattering geometry, its azimuthal integral evaluates to $2\pi$. 

For a specific macroscopic interaction channel $\chan = (l_1, l_2, l_3)$, we define the evaluated innermost angular block $\Omega_{scat}^{(\chan)}(\cmenergy, \fracrel, \halfangle)$ as:
\begin{equation}
\begin{aligned}
\Omega_{scat}^{(\chan)}(\cmenergy,\fracrel,\halfangle) &= \int_{-1}^{1} d(\cos\defangle) B(\urel, \cos\defangle) \Bigg[ \left( \int_{0}^{2\pi} d\aziangle \radbasis_{k_1, l_1}(\vmag') P_{gain}^{(\chan)}(\vpostunit, \halfangle) \right) - 2\pi \radbasis_{k_1, l_1}(\vmag) P_{loss}^{(\chan)}(\halfangle) \Bigg]
\end{aligned}
\end{equation}
where $P_{gain}^{(\chan)}$ and $P_{loss}^{(\chan)}$ represent the geometric Wigner-Eckart filters.

\subsection{Sum-Factorized Duffy Patches}

With the inner scattering manifold defined, the evaluation of the reduced physical tensor $\Rtens$ is expanded as the sum of the two transformed Duffy patches. The Fubini cascade isolates the variables for each patch:

\textbf{Patch 1 Evaluation ($\halfangle = \fracrel \cdot \duffy$):} In the lower triangle, the radial physics depend on the outer coordinate $\fracrel$, allowing the pre-collision state evaluations to be factored outside the inner $\duffy$-loop. The tensor contribution evaluates as:
\begin{equation}
\Rtens_{patch1}^{(\chan)} = \int_{0}^{\infty} d\cmenergy \cmenergy^{\vhs/2}e^{-\cmenergy} \int_{0}^{1} d\fracrel \mathcal{K}_{rad}^{(\chan)}(\cmenergy,\fracrel) \left[ \int_{0}^{1} d\duffy \left( \frac{\fracrel}{\halfangle} \right) \Omega_{scat}^{(\chan)}(\cmenergy, \fracrel, \fracrel \cdot \duffy) \right]
\end{equation}

\textbf{Patch 2 Evaluation ($\fracrel = \halfangle \cdot \duffy$):} In the upper triangle, the angular Wigner-Eckart filters embedded within $\Omega_{scat}$ depend on the outer coordinate $\halfangle$. Because $\fracrel$ changes during the inner loop, the radial polynomials must be evaluated inside the $\duffy$-loop:
\begin{equation}
\Rtens_{patch2}^{(\chan)} = \int_{0}^{\infty} d\cmenergy \cmenergy^{\vhs/2}e^{-\cmenergy} \int_{0}^{1} d\halfangle \left( \frac{\halfangle}{\fracrel} \right) \left[ \int_{0}^{1} d\duffy \Omega_{scat}^{(\chan)}(\cmenergy, \halfangle \cdot \duffy, \halfangle) \mathcal{K}_{rad}^{(\chan)}(\cmenergy, \halfangle \cdot \duffy) \right]
\end{equation}
where the radial weight manifold $\mathcal{K}_{rad}^{(\chan)}$ absorbs the reconstructed pre-collision states and the singularity remainders derived from the Cartesian measure:
\begin{equation}
\mathcal{K}_{rad}^{(\chan)}(\cmenergy,\fracrel) = 2^{(\vhs/2-3)} \fracrel^\vhs \cmenergy^2 (1-\fracrel^2)^2 e^{-\cmenergy\fracrel^2} \radbasis_{k_2, l_2}(\vmag) \radbasis_{k_3, l_3}(\wmag)
\end{equation}

\subsection{Precomputation of Angular Manifolds}

The radial polynomials evaluate via contiguous scalar arithmetic; however, the geometric Wigner-Eckart cross-sections involve trigonometric evaluations. To prevent redundant calculations and ensure contiguous memory access during the tensor contractions, we precompute the angular scattering manifolds on the discrete grid before the execution phase begins. 

These static geometric arrays are constructed as:
\begin{itemize}
    \item \textbf{Patch 1 Geometric Surfaces:} Because the physical incidence angle $\halfangle = \fracrel \cdot \duffy$ changes with both integration variables, the azimuthally-integrated gain and loss filters are precomputed across the 2D sub-domain. To align with column-major memory architecture, these are stored as arrays of size $[N_{\duffy 1}, N_{\fracrel 1}, N_Q]$, ensuring the inner loop variable $\duffy$ acts as the fastest-moving index.
    \item \textbf{Patch 2 Geometric Vectors:} Because the incidence angle $\halfangle$ governs the outer integration loop, it is independent of the inner $\duffy$-loop. The 2D angular surface collapses into a 1D geometric vector of size $[N_{\halfangle 2}, N_Q]$, evaluated at the outer nodes $\halfangle_i$, reducing the memory bandwidth required during the inner contraction.
\end{itemize}

\subsection{Block-Decoupled Execution Sequence}

With the static arrays initialized, the 5D integration is evaluated. The outermost integration over the center-of-mass energy $\cmenergy$ is decoupled, allowing for parallel block evaluation. For a fixed energy node $\cmenergy_i$, the independent patches are contracted as follows:

\begin{enumerate}
    \item \textbf{Patch 1 Contraction (Outer $\fracrel$, Inner $\duffy$):} The algorithm iterates over the fractional coordinate $\fracrel$. Inside this loop, the radial polynomials remain constant. The integration executes the inner $\duffy$-loop, aggregating the 3D angular physics into a dense 1D vector. Exiting the $\duffy$-loop, this intermediate vector is contracted against the static radial polynomials via dense matrix-matrix multiplication (GEMM), accumulating into the tensor $\Rtens$.
    \item \textbf{Patch 2 Contraction (Outer $\halfangle$, Inner $\duffy$):} The algorithm iterates over the angular coordinate $\halfangle$. Here, the geometric Wigner-Eckart filters remain constant. Inside the inner $\duffy$-loop, the pre-collision kinematics $\vmag$ and $\wmag$ are reconstructed dynamically, and the collision kernel is evaluated to reduce the radial physics into a 1D vector. Exiting the $\duffy$-loop, this sum-factorized radial vector is contracted against the static 1D scattering geometry via matrix-vector multiplication (GEMV), preserving the minimal arithmetic sequence.
\end{enumerate}

\section{Derivation of Quadrature Grid Exactness Bounds} \label{app:quad_bounds}
To maintain the spectral convergence of the Galerkin framework, the multi-dimensional quadrature grid must resolve the highest polynomial degree generated by the five-dimensional integrand. The basis functions are constructed from associated Laguerre polynomials and solid spherical harmonics. Given the spectral truncation limits $\kmax$ and $\lmax$, a single particle state represents a polynomial of maximum degree $D_{state} = 2\kmax + \lmax$ in the Cartesian velocity components. 

The total Boltzmann collision integral evaluates the product of three such states (test, target, and incident) alongside the transformed Cartesian measure. Because the 2D Duffy transformation couples the kinematic variables, we derive the resulting maximum algebraic degree along each integration axis.

\subsection{Coordinate-Wise Degree Derivation}

\textbf{Energy Coordinate ($\cmenergy$):} The kinematic rotation dictates that the scalar energies $\vmag^2$ and $\wmag^2$ are linear with respect to $\cmenergy$. Tracing the dependencies through the three interacting states and the integration measure yields a total algebraic degree in $\cmenergy$ of $3\kmax + \lfloor 1.5\lmax \rfloor + 2$. A Generalized Gauss-Laguerre rule absorbs the analytic $\cmenergy^{\vhs/2}e^{-\cmenergy}$ factor. To integrate the remaining polynomial of degree $D$, a Gauss-type rule requires $2N - 1 \ge D$. Therefore, the baseline resolution is:
\begin{equation}
    N_{\cmenergy} = \left\lceil \frac{3\kmax + 1.5\lmax + 3}{2} \right\rceil.
\end{equation}

\textbf{Deflection and Azimuthal Angles ($\defangle, \aziangle$):} The rotational scattering physics remain independent of the radial singularity mapping. Only the post-collision target state depends on the polar deflection, requiring a Gauss-Legendre rule sized at:
\begin{equation}
    N_{\cos\defangle} = \kmax + \lceil 0.5\lmax \rceil + 1.
\end{equation}
The azimuthal angle $\aziangle$ resolves the highest Fourier mode. For the periodic domain $[0, 2\pi]$, the Trapezoidal rule functions identically to a Gauss rule, providing exactness with:
\begin{equation}
    N_{\aziangle} = 2\kmax + \lmax + 1.
\end{equation}

\subsection{Sum-Factorized Duffy Patches}
Before the Duffy split, we determine the independent Cartesian degrees of the fractional coordinate $\fracrel$ and the incidence half-angle $\halfangle$:
\begin{itemize}
    \item \textbf{Degree in $\fracrel$ ($D_{\fracrel}$):} The three interacting states and the transformed Cartesian measure $(1-\fracrel^2)^2$ contribute a maximum degree of $D_{\fracrel} = 6\kmax + 3\lmax + 4$.
    \item \textbf{Degree in $\halfangle$ ($D_{\halfangle}$):} The test state, the geometric test filter $P(1-2\halfangle^2)$, and the $4\halfangle \, d\halfangle$ Jacobian contribute a maximum degree of $D_{\halfangle} = 2\kmax + 3\lmax + 1$.
\end{itemize}

Because the mapped variables decouple the integration loops, we tailor the exactness bounds for the topology of each patch:

\textbf{Patch 1 ($\fracrel > \halfangle$):} We map $\halfangle = \fracrel \cdot \duffy$. The integration measure transforms with an additional Jacobian factor of $\fracrel$ ($d\fracrel \, d\halfangle \to \fracrel \, d\fracrel \, d\duffy$). The highest generic cross-term $\fracrel^{D_{\fracrel}} \halfangle^{D_{\halfangle}}$ maps to $\fracrel^{D_{\fracrel} + D_{\halfangle} + 1} \duffy^{D_{\halfangle}}$. 
\begin{itemize}
    \item \textbf{Outer Loop ($\fracrel$):} Absorbs both variables, requiring exactness for $D_{\fracrel} + D_{\halfangle} + 1 = 8\kmax + 6\lmax + 6$. Solving $2N - 1 \ge 8\kmax + 6\lmax + 6$ yields the outer resolution:
    \begin{equation}
        N_{\fracrel 1} = 4\kmax + 3\lmax + 4.
    \end{equation}
    \item \textbf{Inner Loop ($\duffy$):} Evaluates the mapped $\halfangle$-dependencies, requiring exactness for $D_{\halfangle} = 2\kmax + 3\lmax + 1$. Solving $2N - 1 \ge 2\kmax + 3\lmax + 1$ yields the inner resolution:
    \begin{equation}
        N_{\duffy 1} = \kmax + \lceil 1.5\lmax \rceil + 1.
    \end{equation}
\end{itemize}

\textbf{Patch 2 ($\halfangle > \fracrel$):} We map $\fracrel = \halfangle \cdot \duffy$. The integration measure transforms with an additional Jacobian factor of $\halfangle$ ($d\fracrel \, d\halfangle \to \halfangle \, d\halfangle \, d\duffy$). The cross-term $\fracrel^{D_{\fracrel}} \halfangle^{D_{\halfangle}}$ maps to $\halfangle^{D_{\fracrel} + D_{\halfangle} + 1} \duffy^{D_{\fracrel}}$.
\begin{itemize}
    \item \textbf{Outer Loop ($\halfangle$):} Absorbs both variables, requiring identical exactness for $D_{\fracrel} + D_{\halfangle} + 1 = 8\kmax + 6\lmax + 6$.
    \begin{equation}
        N_{\halfangle 2} = 4\kmax + 3\lmax + 4.
    \end{equation}
    \item \textbf{Inner Loop ($\duffy$):} Evaluates the mapped $\fracrel$-dependencies, requiring exactness for $D_{\fracrel} = 6\kmax + 3\lmax + 4$. Solving $2N - 1 \ge 6\kmax + 3\lmax + 4$ yields:
    \begin{equation}
        N_{\duffy 2} = 3\kmax + \lfloor 1.5\lmax \rfloor + 3.
    \end{equation}
\end{itemize}

\subsection{Selective Padding for Non-Polynomial Components}
While the bounds derived above guarantee integration of the algebraic components, the physical integrand contains two non-polynomial factors: the residual Gaussian envelope $e^{-\cmenergy \fracrel^2}$ and the fractional scattering kernel $\urel^{\vhs}$. Because the Duffy transformation ensures these terms are $C^\infty$ smooth and analytic on the mapped domains, Gauss-Legendre quadrature provides exponential spectral convergence. To resolve the infinite Taylor series of these terms down to machine precision, targeted padding parameters must be added to the theoretical baselines.

This padding is not applied globally across all five dimensions. The rotational scattering angles ($\cos\defangle$ and $\aziangle$) evaluate polynomial or periodic functions. Because they are isolated from the radial singularity mapping, their baseline bounds ($N_{\cos\defangle}$ and $N_{\aziangle}$) provide exact integration without padding. 

Padding is reserved for the variables coupled by the center-of-mass rotation and the 2D Duffy transformation:
\begin{itemize}
    \item \textbf{Energy Coordinate ($\cmenergy$):} The rotation to center-of-mass coordinates leaves behind the residual Gaussian $e^{-\cmenergy \fracrel^2}$. Resolving this exponential decay requires applying the global radial padding parameter ($N_{rad\_pad}$) directly to $N_{\cmenergy}$.
    \item \textbf{Patch 1 ($\fracrel > \halfangle$):} The residual Gaussian and fractional square root dictate the physics. The outer loop integrates the primary radial fractional coordinate $\fracrel$, requiring the radial padding ($N_{\fracrel 1} + N_{rad\_pad}$). The inner $\duffy$ loop evaluates the mapped half-angle dependencies, taking the angular padding ($N_{\duffy 1} + N_{ang\_pad}$).
    \item \textbf{Patch 2 ($\halfangle > \fracrel$):} The physical roles invert. The outer loop integrates the incidence half-angle $\halfangle$, requiring the angular padding ($N_{\halfangle 2} + N_{ang\_pad}$). The inner $\duffy$ loop evaluates the mapped radial fractional dependencies, requiring the radial padding ($N_{\duffy 2} + N_{rad\_pad}$).
\end{itemize}